\def\ZZ         {{\mathbb Z}}

\def\RR         {{\mathbb R}}
\def\CC         {{\mathbb C}}
\def\QQ         {{\mathbb Q}}
\def\PP         {{\mathbb P}}
\def\TT         {{\mathbb T}}
\def\NN         {{\mathbb N}}
\def\A          {{\cal A}}
\def\X          {{\cal X}}
\def\D           {{\cal D}}
\def\J          {{\cal J}}

\def\L          {{\cal L}}
\def\M          {{\cal M}}
\def\O          {{\cal O}}
\def\U           {{\cal U}}
\def\S           {{\cal S}}  
\def\V           {{\cal V}}

\def\rk         {{\rm rk}}

\def\dim          {{\rm dim}}
\def\log         {{\rm log}}

\def\Char       {{\rm Char}}
\def\Spec       {{\rm Spec}}
\def\Ch         {{\rm Ch}}

\def\dim        {{\rm dim}}
\def\log        {{\rm log}}

\def\Ker        {{\rm Ker}}
\def\Im         {{\rm Im}}

\def\cal        {\mathcal}

\documentclass{amsart}
\usepackage{amssymb}
\usepackage{verbatim}
\usepackage{eucal}
\usepackage{mathrsfs}
\usepackage{graphicx}
\usepackage{psfrag}
\usepackage{amsthm}

\newtheorem{theorem}{Theorem}[section]
\newtheorem{lemma}[theorem]{Lemma}
\newtheorem{prop}[theorem]{Proposition}
\newtheorem{corollary}[theorem]{Corollary}

\theoremstyle{definition}
\newtheorem{dfn}[theorem]{Definition}

\theoremstyle{remark}
\newtheorem{remark}[theorem]{Remark}

\begin{document}

\title{Non vanishing loci of Hodge numbers of local systems.}

\author{Anatoly Libgober}
\footnote{The author supported by National Science Foundation grant.}
\address{Department of Mathematics\\
University of Illinois\\
Chicago, IL 60607}
\email{libgober@math.uic.edu}

\begin{abstract} We show that closures of families of unitary local systems
on quasiprojective varieties for which the dimension of a graded
component of Hodge filtration has a constant value can be identified with 
a finite union of polytopes. We also present a local version of this theorem. 
This yields the ``Hodge decomposition'' 
of the set of unitary local systems with a non-vanishing 
cohomology extending Hodge decomposition 
of characteristic varieties of links of plane curves studied 
by the author earlier. 
We consider a twisted version of the characteristic varieties 
generalizing the twisted Alexander polynomials.
Several explicit calculations for complements to 
arrangements are made.
\end{abstract}

\maketitle

\section{Introduction}

Let $\rho: \pi_1(X) \rightarrow U_N$
be a unitary $N$-dimensional representation 
of the fundamental group of a non-singular quasiprojective 
variety $X$ or, in other words, a unitary rank $N$ local system.
It is well 
known (cf.\cite{Timm},\cite{Arapura})
that with this data one can associate 
the cohomology groups $H^i(X,\rho)$  and that 
the cohomology of a unitary local system $V_{\rho}$ supports
canonical mixed Hodge structure, part of which is the Hodge filtration
$...F^{p+1}H^n(V_{\rho})  
\subseteq F^pH^n(V_{\rho}) \subseteq ...$.

One of the purposes of this paper is to study the structure of the 
family $S^{n,p}_{\rho ,l}$ of local systems $V_{\rho} \otimes L_{\chi}$
with a fixed $\rho$ corresponding to unitary characters $\chi: 
\pi_1(X) \rightarrow \CC^*$ 
for which the dimension of graded components
$Gr_F^p=F^p/F^{p+1}$ of the graded space associated 
with the Hodge filtration on cohomology satisfies 
$\dim Gr_F^pH^n(X,V_{\rho \otimes \chi}) \ge l$.
This structure, it turns out, is rather different 
in projective and quasiprojective cases and the difference
will be explained shortly in this introduction.

Secondly,  we derive a {\it local} version of results on the structure
of the sets of local systems with fixed dimension of 
cohomology and also the sets of local systems with 
fixed dimension of associated graded for Hodge filtration groups
(a construction of a mixed Hodge structure in this case, 
apparently absent in the literature and therefore we give one in 
 section \ref{localtranslates}).  
The local and global situations are closely related and 
in previous paper \cite{ample} we showed how such type of local data 
can be used to derive information about sets of local systems 
with fixed dimension of cohomology in some quasiprojective cases.  
For simplicity, we assume in this paper that there are no non-trivial 
rank one local systems on a non-singular compactification $\bar X$ of $X$
i.e. $H^1(\bar X,\CC^*)=0$.

Let us first describe our local results, (which global 
counterparts are in \cite{GL}, \cite{Simpson} and \cite{Arapura})
and consider a pair $(\X,\D)$ consisting of a 
germ of a complex space $\cal X$ having as its  
link a simply connected manifold and a divisor 
${\cal D}=\bigcup D_i$ on $\X$ which is a union of $r$ irreducible components.
In this case the family of the rank one local 
systems on $\X-\D$ 
is  parameterized by a $r$-dimensional affine torus:
$H^1({\cal X}-{\cal D},\CC^*)={\CC^*}^r$ (cf. \cite{DL}).
In section \ref{localtranslates} we show the following:
\begin{theorem}\label{translates}
The collection of local systems 
$$S_l^n=\{\chi \in  H^1({\cal X}-{\cal D},\CC^*) 
\vert \dim H^n({\cal X}-{\cal D},{L}_{\chi}) \ge l\}$$
is a finite union of translated by points of finite order 
sub-tori of $H^1({\cal X}-{\cal D},\CC^*)$.
\footnote{The assumption of simply connectedness of the 
link of $\cal X$ is used only to 
simplify the exposition and can be dropped after some modification 
in the statement. Cf. work 
\cite{budur} for such modification in related context.} 
\end{theorem}

This is a local counterpart of the translated subgroup 
property of the families of rank one local systems 
with jumping cohomology on quasiprojective manifolds (cf. \cite{Arapura}).
Note that ${\cal X}-{\cal D}$ is 
Stein but not quasiprojective and so, except for occasional cases 
when it is homotopy equivalent to a quasiprojective manifold,
this does not follow immediately from \cite{Arapura}. Rather, we
use a Mayer-Vietoris spectral sequence for the union 
of tori bundles on quasiprojective manifolds, which, as we show, 
degenerates in term $E_2$ and this 
together with \cite{Arapura} yields the claim.
As in \cite{innc} (cf. also \cite{DM}), we call $S_l^n$ the 
characteristic variety of $({\cal X},{\cal D})$.

The theorem \ref{translates} was conjectured in \cite{innc} where 
we also conjectured a procedure for calculation of the translated sub-tori 
from the theorem \ref{translates} 
in the case when $\cal X$ is non-singular and 
$\cal D$ is an isolated non-normal crossing divisor 
(see section \ref{inncsection}
for discussion of this case).
The idea in \cite{innc} was to calculate the 
dimensions $h^{p,q,n}_{\chi}$ of the vector spaces defined 
in terms of the mixed Hodge structure on 
the cohomology of a local system $L_{\chi}$, 
using their relation to the homology 
of abelian covers and then to apply the results to obtain the 
translated sub-tori in the theorem \ref{translates} as 
the Zariski closure in $H^1({\cal X}-{\cal D},\CC^*)$ 
of the image of the jumping loci of the Hodge group under the 
exponential map.
A special case of this approach is the algebro-geometric
calculation of the zeros of multivariable Alexander polynomials 
carried out in \cite{alexhodge}. 
In section \ref{localtranslates} we address the issue of the existence of 
a Mixed Hodge structure on the cohomology of local systems
on ${\cal X}-{\cal D}$,
which admit logarithmic extension with eigenvalues belonging to $[0,1)$.
This is a local version of the result of Timmescheidt \cite{Timm} 
(also \cite{Arapura})
and in the case of trivial local systems reduces to 
construction in \cite{DH}. The Hodge numbers of local systems 
corresponding to $\chi$'s having finite order in $H^1(\X-\D,\CC^*)$
obtained from the construction in section \ref{localtranslates}  
coincide with the Hodge numbers $h^{p,q,n}_{\chi}$ 
which are the dimensions of $\chi$-eigenspaces 
of the covering group acting on the Deligne's Hodge groups
of finite abelian covers of $\X-\D$.
This allows us in section \ref{inncsection} to use
the information obtained in section 
\ref{localtranslates} to study the  
Hodge decomposition of abelian covers corresponding to isolated
non normal crossings.

Now let us explain the difference between the loci 
with fixed dimension of $Gr_F^pH^n$ in projective 
and quasiprojective cases. 
On compact Kahler manifolds, the families of local systems 
can be identified with the families of topologically trivial 
\footnote{one can also work with bundles for which the 
first Chern class is a torsion}
holomorphic line bundles (cf. \cite{Simpson} or \cite{budur})
and those with jumping Hodge numbers are the unions
of translated abelian subvarieties of the Picard groups (cf. \cite{GL}). 
As an example in quasiprojective case, let us 
look at the cohomology of local systems on $\PP^1-(0,1,\infty)$.
The unitary local systems are parameterized by the maximal compact 
subgroup in $\CC^* \times \CC^*$. One can show that 
the cohomology $H^1(\PP^1-(0,1,\infty),{L}_{\chi})$ 
of a local system corresponding to 
$\chi: \pi_1(\PP^1-(0,1,\infty))
\rightarrow \CC^*$ and having a finite order $n$ is isomorphic to 
the $\chi$-eigenspace $H^1(X_{n,n},\CC)_{\chi} \subset H^1(X_{n,n},\CC)$ 
where $X_{n,n}$ is the abelian cover of $\PP^1-(0,1,\infty)$ corresponding to 
the reduction modulo $n$: 
$$H_1(\PP^1-(0,1,\infty),\ZZ)=\ZZ \oplus \ZZ
\rightarrow \ZZ/n\ZZ \oplus \ZZ/n\ZZ$$
(in other words: $H^1(X_{n,n},\CC)_{\chi}=\{ a \in H^1(X_{n,n},\CC)
 \vert g \cdot a=\chi(g)a \ {\rm where}\  g \in \ZZ/n\ZZ \oplus \ZZ/n\ZZ\}$). 
Such an abelian cover has as a model
the complement to the fixed points of non identity elements in the 
group $\ZZ/n\ZZ \oplus \ZZ/n\ZZ$ acting on a Fermat curve: 
$$F_n: \ x^n+y^n=z^n, \ \ \ (x,y,z) \rightarrow (\zeta_n^ax,\zeta_n^by,z)$$
($\zeta_n$ is a primitive root of degree $n$, $a,b \in \ZZ/n\ZZ$).
In particular, each $H^1(\PP^1-(0,1,\infty),{L}_{\chi})$
acquires the Hodge and weight 
filtrations induced from $H^1(X_{n,n},\CC)_{\chi}$.
We shall focus on the weight 1 
component $Gr^W_1H^1(X_{n,n})=H^1(F_n)$ and 
will look for those $\chi$ for which 
$$Gr^1_FGr^W_1H^1(\PP^1-(0,1,\infty),L_{\chi})=
Gr^1_FGr^W_1H^1(F_n,\CC)_{\chi} \ne 0$$
Now the cohomology classes $H^{1,0}(F_n,\CC)$ are 
represented by the residues of meromorphic 2-forms on $\PP^2$
with poles {\it having order one} along $F_n$ {cf. \cite{griffiths}:
$${P(x,y,z) zdx \wedge dy} \over {x^n+y^n-z^n}$$
where $P(x,y,z)$ is a homogeneous polynomial of degree $n-3$. 
Since the group $\ZZ/n\ZZ \oplus \ZZ/n\ZZ$ 
acts via multiplication of the coordinates (in particular generators
act as follows:   
$g_x: (x,y,z) \rightarrow (\zeta_nx,y,z), g_y: (x,y,z) \rightarrow
(x,\zeta_ny,z)$) each
non zero eigenspace of the character $\chi_{a,b}$ such that 
$\chi_{a,b}(g_a)=e^{{2 \pi i (a+1)}\over n},
\chi_{a,b}(g_b)=e^{{2 \pi i (b+1)} \over n}$ 
is generated by the monomial forms
${{x^ay^bz^{n-3-a-b}z dx \wedge dy} \over
 {x^n+y^n-z^n}}
 \ \ (0 \le a+b \le n-3)$.
On the universal cover of the space of unitary characters, they are
represented in terms of $\bar a=a+1, \bar b=b+1$
by the triangle $$\{{\bar a \over n}, {\bar b \over n} \vert \bar a > 0 \ 
\bar b > 0, \bar a+\bar b  < n \}
\footnote{One can check that weight two part 
$Gr^1_FGr^W_2H^1(X_{n,n})_{\chi}$ adds
characters on the boundary 
of this triangle.}$$ 
 
The situation described in this example is quite general. 
In section \ref{quasiproj} we prove the following:

\begin{theorem}\label{quasiprojhodge} Let $X$ be a quasiprojective manifold 
without non-trivial rank one local systems on 
a non singular compactification, 
$\rho: \pi_1(X) \rightarrow U_N$ be a $N$-dimensional unitary 
representation , $\Char \pi_1(X)$ 
be the torus of characters of the fundamental group and $\Char^u\pi_1(X)$ 
be the subgroup of unitary characters. Let $\U$ be the fundamental 
domain of $\pi_1(\Char^u \pi_1(X))$ acting on the universal cover 
$\widetilde {\Char^u \pi_1(X)}$ of the torus $\Char^u \pi_1(X)$
and $exp: \U \rightarrow \Char^u \pi_1(X)$ be 
the universal covering map
\footnote{We identify the universal cover with the tangent space
to $\Char^u\pi_1(X)$, the universal covering map with the 
exponential map and the fundamental domain $\cal U$ with the unit cube in
the tangent space.}. 
 Then 
\begin{equation}
\S^{n,p}_{\rho,l}=
\{\chi \in \Char^u \pi_1(X) \vert \dim 
Gr_F^pH^{n}({V_{\rho} \otimes L_{\chi}}) \ge l\}
\end{equation}
is a finite union of polytopes in $\U$ i.e. the subsets, each of which 
is the set 
of solutions for a finite set of inequalities $L \ge 0$ where $L$ is a 
linear function.
\end{theorem}

The argument is based on a study of Deligne extensions and yields
 an independent proof of quasiprojective version of theorem 
\ref{translates}. (i.e. the result of \cite{Arapura}).

We also prove that in the local situation of theorem \ref{translates}
the structure of the families local system with fixed 
Hodge numbers is similar to the one in 
quasiprojective case as described in theorem \ref{quasiprojhodge}.
More precisely, we show the following:

\begin{theorem}\label{Main} Let $(\X,\D)$ be a germ 
of a divisor with $r$ irreducible components on a germ 
of a complex space as in theorem \ref{translates}. 
Let $X=\X-\D$  and $Char X=Char \pi_1(X)$ be the space of
characters of the fundamental group 
of $X$. 
\par \noindent 
Let $Char^u \pi_1(X) \subset Char \pi_1(X)$ be the maximal torus of 
unitary characters and let ${\cal U} \subset 
\widetilde {\Char^u \pi_1(X)}$
 be a fundamental domain of the covering 
group of the universal cover $\widetilde {\Char^u \pi_1(X)}$
of $Char^u \pi_1(X)$.
For each pair $(p,n)$ and $l \in \NN$ there exists a collection 
of polytopes
$\S_{l}^{p,n} \subset {\cal U}$ such that 
the image of $\S_{l}^{p,n}$  
under the covering map ${\cal U} \rightarrow Char^u\pi_1(X)$
consists of the local systems $L$ satisfying $rk F^p/F^{p+1}H^n(L)=l$.
The image of the union 
$\cup_{p,(...,l_i,...),\sum l_p=l} \S_{l_i}^{n,p}$ under the covering map 
$\widetilde {\Char^u \pi_1(X)} \rightarrow Char^uX$ 
is the unitary part of characteristic variety: $S_l^n \cap \Char^u\pi_1(X)$.
\end{theorem}

One of the tools used in the proof of theorem \ref{Main} is 
the Mayer Vietoris spectral sequence for the cohomology 
of rank one local systems on the union of quasiprojective 
manifolds with normal crossings which is discussed in lemma
\ref{degenerationMV}. 

In the case when $(\X,\D)$ is an isolated non-normal 
crossing (which is studied in \cite{innc} and \cite{DL})
this has as a consequence the fact that the support of the homotopy 
groups of the complements has canonical decomposition into a union
of polytopes. These polytopes can be related to the distribution
of the Hodge numbers of abelian covers and will be discussed
in section \ref{inncsection}. 

The next section studies an abelian 
generalization of the twisted Alexander polynomials studied in \cite{CF}.
We also prove a generalisation of the cyclotomic property of the 
roots of Alexander polynomials (cf. theorem \ref{roots}) 
which lead to a restrictions 
on class of groups which are isomorphic to the fundamental groups
of the complements to plane algebraic curves (cf. theorem \ref{restriction}).
Finally, in the last section,
 we make explicit calculations 
of  polytopes $\S^{n,p}_l$ in several examples. 

The author wants to thank N.Budur and J.I.Cogolludo for their
comments on a draft of this paper.
Note that the polytopes \footnote{called ``polytopes of quasiadjunction''}
 in this circle of questions were introduced in \cite{abcov}
(cf. also \cite{alexhodge}). Polytopes 
in the quasiprojective case 
were studied recently by Nero Budur in \cite{budur} 
using Mochizuki's work. 

\section{Hodge numbers of local systems on quasiprojective varieties}
\label{quasiproj}
Let $X$ be a quasiprojective manifold and 
let $\bar X$ be a compact projective manifold such that 
$H^1(\bar X,\CC^*)=0$, $\bar X-D=X$ where $D$ is a divisor  
with normal crossings.
Denote by $\Omega^1_{\bar X}({\rm log} D)$ the sheaf 
of logarithmic 1-forms.
We shall fix a $N$-dimensional unitary representation
$\rho: \pi_1(X) \rightarrow U_N$  
such that there is a locally trivial bundle $\cal V$ on $\bar X$,
a meromorphic connection 
$\nabla: \V \rightarrow \Omega^1_{\bar X}({\rm log} D)
\otimes \V$ for which the restriction on $X$ is flat and the 
corresponding holonomy representation is $\rho$ (cf. \cite{delignediff}).
Let $\chi: \pi_1(X) \rightarrow \CC^*$ be a unitary character of the 
fundamental group. Denote by  $V_{\rho \otimes \chi}$  the
local system corresponding to representation $\rho \otimes \chi$. 
It follows from 
\cite{Timm} and \cite{Arapura} that the cohomology groups 
$H^i(X,V_{\rho \otimes \chi})$ support 
three filtrations $F,\hat F, W$ ($F,\hat F$ being decreasing and 
$W$ increasing) 
such that 
$$\dim Gr^W_nH^i(X,V_{\rho} \otimes L_{\chi})
=\dim F^pGr_n^W + \dim \hat F^qGr_n^W \ \ \ (p+q=n)$$
Moreover, these filtrations on the cohomology are independent of 
a particular compactification $\bar X$ used in their construction
but rather depend only on $X$
and $\rho$.
Note that $F$-filtration is the one resulting 
from the degenerating in $E_1$-term Hodge-deRham spectral sequence
(cf. \cite{Timm}):
\begin{equation}\label{hodgespectral}
E_1^{p,q}=H^p(\Omega^q(\log D)\otimes \V_{\rho \otimes \chi}) 
\Rightarrow H^{p+q}(V_{\rho \otimes \chi}).
\end{equation}
Here $\V_{\rho \otimes \chi}$ 
is a vector bundle on $\bar X$ with flat logarithmic connection
$\nabla$ such that the eigenvalues of the residues matrices 
$Res \nabla$ along components of $D$ 
satisfy $0 \le {\rm Re} <1$
and such that on $X$ the holonomy of $\nabla$ is  the representation 
$\rho \otimes \chi$. We call such $\V$ the Deligne's (canonical) extension
of a bundle on $X$ supporting the local system.

The set of unitary rank one local systems $\Char^u(\pi_1(X))$ 
is parameterized by the maximal 
compact subgroup of the torus $Hom(\pi_1(X),\CC^*)$. We assume
for simplicity {\it throughout this paper} 
that this torus is connected (i.e., $H_1(X,\ZZ)$ is 
torsion free) leaving to an interested reader the general case.
As already was mentioned, we view the 
universal cover $\widetilde {\Char^u(\pi_1(X))}$ 
as the tangent space to $\Char^u(\pi_1(X))$
at the identity 
with the covering map  
being the exponential map. A choice of generators in $H_1(X,\ZZ)
=\oplus \ZZ^r$ allows one to identify the universal 
cover $\widetilde {\Char^u(\pi_1(X))}$
with $\RR^n$ so that  the generator of 
the i-th summand acts as $(x_1,...,x_r) \rightarrow (x_1,..,x_i+1,...
x_r)$. An exponential map yields the lattice 
$\widetilde {\Char^u(\pi_1(X))}_{\ZZ}=\Ker( exp ) \subset 
\widetilde {\Char^u(\pi_1(X))}$.
We use 
the unit cube 
\begin{equation}\label{funddomain}
\U: \ \ \  \{(x_1,...,x_r) \vert 0 \le x_i <1\}
\end{equation}
as the fundamental domain for the action of 
$\widetilde {\Char^u(\pi_1(X))}_{\ZZ}=H_1(X,\ZZ)$ 
on the universal cover 
$\widetilde {\Char^u(\pi_1(X))}$.
A subset of $\U$ is called {\it a polytope} if 
there exists a collection of linear functions $l_j \in Hom (
\widetilde {\Char^u(\pi_1(X))}_{\ZZ},\QQ)$ 
on the universal cover $\widetilde {\Char^u(\pi_1(X))}$ 
so that this subset is the set of 
solutions of a finite collection of inequalities $a_j \le  l_j(u) < a_j' \ \  
(a_j,a_j' \in \RR)$.

The main result in this section is the following:

\begin{theorem}\label{quasiproj}
 Let $X$ be a quasiprojective variety and let $\bar X=X \cup D$ 
be a compactification 
such that $H^1(\bar X,\CC^*)=0$. 
Let $S^{n,p}_{\rho,l}$ 
be the subset 
of $\U$ defined as follows.
\begin{equation}\label{slocus}
S^{n,p}_{\rho,l}=
\{u \in \U \vert \dim Gr_F^pH^{n}( V_{\rho \otimes {exp(u)}}) \ge l\}
\end{equation}
Then (\ref{slocus}) is a finite union of polytopes.
 \end{theorem}

\begin{proof} Let $\pi: (\tilde X, \tilde D)
 \rightarrow (\bar X, D)$ 
be a log-resolution of $(\bar X,D)$ i.e. the total transform 
$\tilde D=\bigcup_{k=1,...,K} \tilde D_k$
of $D$ is a normal crossing divisor. 
We claim that on $\tilde X$  
there are only finitely many bundles $\V$ which are the  
Deligne's extensions of 
the local systems in the family $V \otimes L_{\chi}$ 
where $\chi \in \Char \pi_1(X)$.
Moreover, the collection 
of unitary local systems in the torus $\Char^u(\pi_1(X))$ 
 having a fixed Deligne's extension 
is the image via the exponential map of a polytope in the 
universal cover  $\widetilde {\Char^u(\pi_1(X))}$. The above theorem 
clearly follows from this.

Let $\vert \vert\omega_{i,j}(\rho) \vert \vert$ (resp. $\omega^{\chi}$)
be the connection matrix of a flat logarithmic 
connection (for a locally trivial bundle $\cal V$ on $\tilde X$)
corresponding to the local system $\rho$ (resp. $\chi$), where 
$\omega_{i,j}({\rho}),
\omega^{\chi} \in \Gamma(U,\Omega^1_{\tilde X}(\log (\tilde D)) $ 
are logarithmic 1-forms in a chart  
$U \subset \tilde X$ biholmorphic to a disk so that $U \cap \tilde D$ 
in $U$ is given by $z_1 \cdot ....\cdot z_k=0$.
For a component $\tilde D_k$ of the normal 
crossing divisor $\tilde D$ let 
$R_{i,j}^{\tilde D_k}(\rho)=Res_{\tilde D_k} \omega_{i,j}(\rho)$ be 
the entries of the matrix of residues along $\tilde D_k$ 
(for the extension $\V$).
Let $\beta^{\tilde D_k}_i(\rho)$ be the collection of eigenvalues of 
the matrix $R_{i,j}^{\tilde D_k}(\rho)$.
The components of the matrix of 
connection corresponding to the local system $V \otimes L_{\chi}$ 
in a
basis $v_i, i=1,...,N$ (resp. $e$) of $V$ (resp. $L_{\chi}$)
are $\omega_{i,j}v_j\otimes e+
\omega^{\chi} v_i \otimes e$ 
(the connection on the tensor 
product is given by 
$\nabla(v \otimes e)=\nabla(v) \otimes e+v \otimes \nabla(e)$).
Hence the eigenvalues of $Res_{D_k}
\nabla_{\rho \otimes \chi}$ 
can be computed as: 
\begin{equation}\label{eigenvals}
\beta^{D_k}_{i}(\rho)+ Res_{\tilde D_k} \nabla_{{\chi}} \ \ \ (i=1,...,N).
\end{equation} 
where $Res_{\tilde D_k} \nabla_{{\chi}}$ is the residue of 
the log-connection $\nabla_{{\chi}}$ in the rank one bundle 
along $\tilde D_k$.

Let $\omega_1,....\omega_r, \ \ \omega_s \in H^0(\tilde X,
\Omega^1_{\tilde X}(\log(\tilde D)))$
be a basis of the space of log 1-forms. Notice that 
$\dim H^0(\tilde X,\Omega^1_{\tilde X}(\log(\tilde D)))=r$ since
$\dim H^1(X,\CC)=\dim H^0(\tilde X,\Omega^1_{\tilde X}(\log(D)))$ $
+\dim H^1(\bar X,\O_{\bar X})$
(cf. \cite{delignehodge}), $H^1(\bar X,\CC^*)=H^1(\tilde X,\CC^*)$ and
we assumed that $H^1(\tilde X,\CC^*)=0$.
We shall select  forms $\omega_s$ so that their cohomology classes 
 all belong to  $H^1(X,\ZZ) 
\subset H^1(X,\CC)$. 
The set of $(\alpha_1,....,\alpha_r) \in \CC^r$ such that 
for a fixed $\rho$ the connection $\nabla^{\alpha_1,...,\alpha_r}_{\rho}$
with the matrix 
\begin{equation}\label{connectionvariation}
\vert \vert \omega_{i,j}^{\alpha_1,....,\alpha_r} \vert \vert=
\vert \vert \omega_{i,j} \vert \vert +(\sum \alpha_s \omega_s)I
\end{equation}
yields a unitary local system form a $r$-dimensional 
$\RR$-subspace $\widetilde {\Char^u(\pi_1(X))}$
of $\CC^r$.
Indeed, all connections with matrix 
$\vert \vert \omega_{i,j}^{\alpha_1,....,\alpha_r} \vert \vert$ 
are flat since the connection corresponding to $\vert \vert \omega_{i,j} \vert \vert $ is flat and 
any holomorphic log-form is closed (cf. \cite{delignehodge}).
Moreover, the holonomy of a connection along a path is given
by solutions of a system of first order ODEs and hence is the exponent 
of  a matrix $A$ depending on the each summand in (\ref{connectionvariation}) 
linearly. Hence those $(\alpha_1,...,\alpha_r)$ which yield a skew-hermitian 
$A$ form a $\RR$-linear subspace in the universal cover 
$\widetilde {\Char (\pi_1(X))}=\CC^r$ of the space of characters. 
In addition, there is a compact connected 
fundamental domain  
$\tilde \U \subset \widetilde {\Char^u(\pi_1(X))}$ for the 
action of $H_1(X,\ZZ)$ 
such that the family of unitary local systems
corresponding  to elements in $\tilde \U$ coincides 
with $\Char^u(\pi_1(X))$. 
One can take $\tilde \U$ so that it 
is an affine transform of $\U$ given by (\ref{funddomain}). 


The eigenvalues of the 
residue matrix of the connection $\nabla_{\rho}^{\alpha_1,...,\alpha_r}$ along 
a component $\tilde D_k$ in an extension $\M$ (i.e. a locally trivial bundle 
on $\tilde X$) are given by:   
  
\begin{equation}\label{sbeta}
 N_i^{\tilde D_k}(\alpha_1,...,\alpha_r,\M)=
 \beta^{\tilde D_k}_i(\rho)+\sum \alpha_s Res_{\tilde D_k}\omega_s
 \ \ \ i=1,...,\rk V
\end{equation}
In particular, for each $i=1,...,rk V$ and $\tilde D_k$,  
the image of the map $\tilde \U \rightarrow \RR$ 
given by  $(\alpha_1,...,\alpha_r ) \rightarrow 
N^{\tilde D_k}_i(\alpha_1,...,\alpha_r )$ is bounded
and for any collection of integers 
$n_i^{\tilde D_k},i=1,...,rk V, k=1,...,K$
 the set $(\alpha_1,...,\alpha_r)$ such that 
\begin{equation}\label{residueinequality}
n_i^{\tilde D_k} \le 
N_i^{\tilde D_k}(\alpha_1,...,\alpha_r) <n_i^{\tilde D_k}+1 
\end{equation} 
is a (possibly empty) polytope in $\tilde \U$.

Now let us describe the bundles on $\tilde X$ which are
 the Deligne's extensions
of the connections in the family given by (\ref{connectionvariation}).
Let us call collection $n_i^{\tilde D_k}$ realizable if there 
is a connection corresponding to $(\alpha_1,...,\alpha_r)$ such that 
$n_i^{\tilde D_k}=[N_i^{\tilde D_k}(\alpha_1,...,\alpha_r)]$. 
Let $\M(n_i^{\tilde D_k})$
be the Deligne extension of this connection. We claim that 
all connections with matrices (\ref{connectionvariation}) and residues 
satisfying (\ref{residueinequality}) have 
 $\M(n_i^{\tilde D_k})$ as its Deligne extension. 
Consider a sufficiently fine cover of $\tilde X$ by open sets 
such that both bundles $\M$ and $\M(n_i^{\tilde D_k})$ can be trivilized. 
Over each open set $U$ 
of such a cover one has the transition matrix $g_U 
\in GL_N(\Gamma(U,\O(*\tilde D))$ from a frame $e_1,...e_N$
of $\M$ to a frame $e_1',,,,e_N'$ 
of $\M({n_i^{\tilde D_k}})$. The connection matrices $\Omega(e_1,..e_N)$
and $\Omega(e_1',..,e_r')$ are related as follows:
\begin{equation}
  \Omega(e_1,..,e_N)=g_U^{-1}dg_U+\Omega(e_1',...e_N')
\end{equation}
Since both $exp(2 \pi i Res_{\tilde D_k} \Omega(e_1,..,e_N))$ 
and $exp(2 \pi i Res \Omega(e_1',..,e_N'))$ conside with the 
monodromy of the connection around $\tilde D_k$ (cf. \cite{delignediff} 
Prop. 3.11)
one sees that the eigenvalues of $Res_{\tilde D_k} g_U^{-1}dg_U$ are integers.
Moreover, these residues are equal to $n_i^{\tilde D_k}$  
since $\M(n_i^{\tilde D_k})$ is the Deligne 
extension and the eigenvalues of the residues of the connection 
$\Omega(e_1',...e_N')$
 are in 
interval $[0,1)$. This also implies that 
in the frame $e_1',...e_N'$ the eigenvalues of the residues of the connections
matrix of any connection satisfying (\ref{residueinequality}) has the residues 
in the interval $[0,1)$ i.e. has $\M(n_i^{\tilde D_k})$ as 
its Deligne extension.

In the case when $N=1$ the relation between the extensions 
$\M$ and $\M(n_i^{\tilde D_k})$ is particularly simple.
Indeed, if $\nabla$ is a logarithmic connection in a bundle $\M$ on 
$\tilde X$ and $B=\sum b_i \tilde D_i$ is a divisor on $\tilde X$ 
having support
on $\tilde X-X$ then $\nabla$ induces a connection $\nabla^B$ on 
$\M \otimes \O_{\tilde X}(B)$. The residues  
$Res(\cdot): \M \rightarrow \M \otimes \O_{\tilde D_j}$ of two connections
related as follows: 
\begin{equation}
Res_{\tilde D_j}(\nabla^B)=
Res_{\tilde D_j}(\nabla)-b_j id \vert_{\tilde D_j}
\end{equation}
(cf. \cite{esnault}, lemma 2.7).
In particular, for each member of the variation of the connection given by 
(\ref{connectionvariation}) the corresponding Delinge's connection 
has the form $\M \otimes \O_{\tilde X}(B)$ where the coefficients $b_j$ 
in $B=\sum b_j\tilde D_j$ chosen so that the residues of $\nabla_{\tilde D_j}$
will satisfy $0\le Res_{\tilde D_j}(\nabla^B) <1$. In the case or 
arbitrary $N$ the above argument shows that  
possible number of isomorphism classes of 
bundles on $\bar X$ which are the Deligne's extensions of the connections
(\ref{connectionvariation}) does not exceed: 
\begin{equation}
\sum_{\tilde D_k,i} 
[sup_{(\alpha_1,...,\alpha_r)}
N^{\tilde D_k}_i(\alpha_1,...,\alpha_r)-inf_{(\alpha_1,...,\alpha_r)}
N^{\tilde D_k}_i
(\alpha_1,...,\alpha_r)]
\end{equation} 
($[..]$ is the integer part).   

From the degeneration of the spectral sequence (\ref{hodgespectral}),
it follows that 
$$Gr^p_FH^k(V_{\rho \otimes \chi})=H^p(\tilde X, \Omega^{k-p}_{\tilde X}
({\rm log} \tilde D) 
\otimes \V_{\rho \otimes \chi} )$$
Since the logarithms of characters $\chi$ yielding a given 
Deligne's extension form a polytope (\ref{residueinequality})
in $\U$ and 
since the rank of $Gr^p_FH^k(V_{\rho \otimes \chi})$ 
depends only on the Deligne's extension (and hence 
has 
only finitely many values bounded by the number of Deligne's extensions 
of the local systems $V_{\rho \otimes \chi}$),
we see that the collection of $\alpha_i$ having 
fixed values is a finite union of polytopes.

It follows from (\ref{sbeta}) and (\ref{residueinequality}) 
that each polytope is given by inequalities:

\begin{equation}\label{Deligneinequalities}
N_i^{\tilde D_k}(\V) \le 
\beta^{\tilde D_k}_i(\rho)+\sum \alpha_s Res_{\tilde D_k}\omega_s
< N_i^{\tilde D_k}(\V)+1
 \ \ \ i=1,...,\rk V, \ \ k=1,...,K
\end{equation}
where $N_i^{\tilde D_k}(\V)$ are integers depending on the Chern classes
of bundle $\V$ on $\tilde X$. 
Since classes of $\omega_s$ belong to $ H^1(X,\ZZ) \subset H^1(X,\CC^*)$
one has $Res_{\tilde D_k}\omega_s \in \QQ$ for any $k$.
It follows that the linear subspaces in $\widetilde{\Char^u(\pi_1(X))}$ 
spanned by the sets of 
solutions of $\rk V \cdot K$ inequalities (\ref{Deligneinequalities})
for fixed $\V$ 
are given by equations with integer relatively prime
coefficients and constant terms having the form $\gamma_1 
\beta^{\tilde D_k}_i +\gamma_2
\ (\gamma_1,\gamma_2 \in \QQ)$.
\end{proof}

\begin{corollary}\label{finiteorder} Let $\rho$ be 
a representation with abelian image and let $\rho_1,....,\rho_N 
\in \Char^u \pi_1(X)$ be its irreducible components. 
The isolated characters  $\chi \in \Char^u(\pi_1(X))$ 
for which ${\rm dim} H^k(X,\rho \otimes {\chi}) \ge l$ generate a subgroup of
$\Char^u\pi_1(X)$ such that its subgroup generated by 
$\rho_i \ (i=1,...,N)$ has finite index. 
\end{corollary}

\begin{proof} Since the system of inequalities (\ref{Deligneinequalities})
has only isolated solutions, the logarithm of an 
isolated character $\chi$ as above, is the solution to 
a finite system of equations. 
As was remarked above,
each linear subspace spanned by the sets of 
solutions is given by equations with integer relatively prime
coefficients and constant terms which can be written as
 $\gamma_1 \beta +\gamma_2,
 (\gamma_1,\gamma_2 \in \QQ)$. 
Such subspace contains a point $(\alpha_1,...,\alpha_r)$
with coordinates having the form $\gamma_1\beta+\gamma_2, \ (\gamma_1, 
\gamma_2 \in \QQ)$.
The exponential map takes such a point in $\widetilde{\Char^u(\pi_1(X))}$
to a point in ${\Char^u(\pi_1(X))}$
with 
coordinates  ${[e^{2 \pi i \beta}]}^{\gamma_1}e^{2 \pi i \gamma_2}$. 
If the representation $\rho$ has an abelian image then 
the subgroup in $\Char^u(\pi_1(X))$ generated 
by characters which are the irreducible component of $\rho$
is the subgroup of the points with the coordinates 
$e^{2 \pi i \beta}$ and the claim follows.
\end{proof}

\begin{corollary}\label{finiteorder2}
The subsets 
$$S^n_{\rho, l}=\{ \chi \in \Char^u \pi_1(X) \vert 
\dim H^n(X,L_{\rho \otimes \chi}) \ge l\}$$  are unions of 
finite collections of 
translated subgroups.
If $\rho$ has an abelian image and 
if $\bar \rho_1,..,\bar \rho_N \in \Char^u(\pi_1(X))$ 
are  the irreducible components of $\rho$, then 
translations can be made 
by points 
generating a subgroup of $\Char^u(\pi_1(X))$ 
containing the subgroup generated by $\rho_1,...,\rho_N$ 
as a subgroup of finite order.
\end{corollary}

\begin{proof} The subsets $S^n_{\rho, l}$ clearly are
algebraic subsets of $\Char^u\pi_1(X)$. On the other hand 
it follows from theorem \ref{quasiproj} that 
$S^n_{\rho,l}=\bigcup_{l_0+...+l_{n}=l} \cap_{p=0,...,n} S^{n,p}_{\rho,l_p}$
is a union of polytopes. The Zariski closure in 
$\Char^u\pi_1(X)$ of the image of 
a polytope   
is the intersection of algebraic subvarieties of  $\Char^u\pi_1(X)$
containing this image. Since a polytope spans a linear subspace in 
 $\widetilde{\Char^u\pi_1(X)}$, the Zariski closure 
is the image of the exponential map restricted on 
the spanning subspace i.e.   
it is a translated subgroup of $\Char^u\pi_1(X)$. The claim about the 
translation points follows by the argument identical with 
the one used in the proof of Corollary \ref{finiteorder}.
\end{proof}

\begin{remark} The proof of the theorem \ref{quasiproj} 
yields an alternative proof of the results in \cite{Arapura}.
In \cite{Arapura} it is shown that (at least for trivial 
$\rho$) similarly defined subgroups of the full group of 
characters of $\pi_1(X)$ are unions of finite collections of 
translated subgroups.
In particular the subgroups $S^n_{l}$ in $\Char^u\pi_1(X)$
determine completely the subgroups in $\Char\pi_1(X)$.
The corollary \ref{finiteorder2} shows that translation is by points 
of finite order (rather than just the unitary ones 
cf. \cite{Arapura}).
This also implies that 
the polytopes $S^{n,p}_{\rho,l}$ from the theorem \ref{quasiproj}
for $\rho=1$
determine the translated subgroups of $\Char \pi_1(X)$
completely. In \cite{ample}, a procedure was outlined 
for calculation of the polytopes corresponding to $p=0$ component 
of the Hodge filtration in some geometrically interesting cases.
 \end{remark}

\begin{remark} Proofs of the corollaries \ref{finiteorder}
and \ref{finiteorder2} use as the essential step the fact that the coordinates
of translation characters belong to cyclotomic extension of the 
field generated by $e^{2 \pi i \beta^{D_k}_i}$ for all $D_k$ and $i$.
This provides information about translation points also in the case 
when ${\rm Im}\rho$ is non-abelian.
\end{remark}

\begin{corollary}\label{restriction}
 Let $i:X_1 \rightarrow X_2$ be an 
embedding of  quasiprojective submanifold $X_1$ into a 
quasi-projective submanifold $X_2$. Let $\rho$ be a
unitary representation of $\pi_1(X_2)$.
Then either 
$$\{\chi \in \Char^u \pi_1(X_2) \vert \dim \Ker H^n(X_2,\L_{\rho \otimes \chi})
\rightarrow H^n(X_1,\L_{i*(\rho \otimes \chi)})=l \} $$
or 
$$\{\chi \in \Char^u \pi_1(X_2) \vert \dim \Im H^n(X_2,\L_{\rho \otimes \chi})
\rightarrow H^n(X_1,\L_{i*(\rho \otimes \chi)})=l \} $$
are the images of a union of polytopes in the universal cover of $U$.
\end{corollary}

\begin{proof} Let us select a resolution of a pair 
$(X_1,X_2)$ (i.e., the compactifying divisor $D_2$ of $X_2$ 
is a normal crossings divisor and the closure of $X_1$ 
in this compactification has at infinity the divisor $D_1$ which 
also is a normal crossings 
divisor).
The Deligne extension for the flat connection 
corresponding to the local system $i^*(\rho \otimes \chi)$ is 
restriction of the Deligne extension of connection of 
$\rho \otimes \chi$. Hence the characters in
the corollary are the characters giving the Deligne extension 
for which the corresponding map $H^q(X_2,\Omega^p(\log D_2) \otimes \V)
\rightarrow H^q(X_1,\Omega^p(\log D_1) \otimes \V)$ has a fixed 
dimension of the kernel (resp. image).
\end{proof}

\section{Local translated subgroup theorem 
for complement to germs of singularities and Hodge 
structure on cohomology of local systems}\label{localtranslates}

In this section we prove the local translated subgroups
theorem \ref{translates}  
and the local counterpart of the theorem \ref{quasiproj}.

We shall use the same notations as in introduction: $\X$ 
is a germ of a complex space with an isolated singularity,
$\D$ is a divisor on $\X$ with arbitrary singularities
having $r$ irreducible components. Assume that 
the link of $\X$, i.e. the intersection with 
a small sphere about the singular point, is simply connected
(cf. \cite{DL}).
This implies that:
 \begin{equation}H_1(\X-\D,\ZZ)=\ZZ^r
\end{equation}
(cf. \cite{innc}, \cite{DL}).
As generators for $H_1(\X-\D,\ZZ)$ one can take 
the classes of 1-cycles each being the boundary
of small disks transversal to $\D$ at a non-singular point
of each component. In particular $H^1(\X-\D,\CC^*)={\CC^*}^r$ 
and hence the rank one local systems are parameterized
by the torus with a fixed coordinate system.

The main results of this section are the following two theorems
(cf. theorem \ref{translates} in the introduction).

\begin{theorem}\label{translatedsbg} 
Let $(\X,\D)$ be a germ of a pair where 
$\X$ has an isolated normal singularity 
with a simply connected link
and $\D$ is  divisor with  
$r$ irreducible components 
$D_i \ (i=1,...r)$.
Let 
$$
{\cal S}^n_{l}=\{\chi \in \Char \pi_1(\X-\D) 
\vert H^n(\X-\D,\L_{\chi}) \ge l \} \ \ (1 \le n \le \dim \X)$$
where $\L_{\chi}$ is the local system corresponding to the character $\chi$.
Then ${\cal S}^n_l$ is a union of a finite collection of 
translated subgroups for any $n$ and $l$. More precisely, 
there are (possibly trivial)
subgroups $T_i \subset \Char \pi_1(\X-\D),\ (i \in {\cal I}, 
{\rm Card}\ {\cal I} < \infty),$ 
and torsion characters
$\rho_i$ such that 
$${\cal S}^n_l=\bigcup_i \rho_iT_i$$ 
\end{theorem}

The arguments we are using in the proof of the theorem 
\ref{translatedsbg}
also yield the existence of the mixed Hodge 
structure on the cohomology of the unitary local systems 
extending the results of \cite{Timm}, \cite{Arapura} 
to the local case:

\begin{theorem}\label{localmhs}  Let $\X$ be a germ 
of an analytic space 
having an isolated normal singularity and let 
$\D$ be a divisor on $\X$. Denote by $\rho$ a {\it unitary} 
representation of $\pi_1(\X-\D)$ and
let $\L_{\rho}$ be the corresponding
local system. Then the cohomology groups $H^i(\X-\D,\L_{\rho})$ support 
the canonical $(\CC)$-mixed Hodge structure compatible with 
the holomorphic maps of pairs $(\X,\D)$ 
endowed with a local system on the complement $\X-\D$.
\end{theorem}

We refer to \cite{Arapura} and \cite{Arapura1} for a discussion 
of $\CC$-mixed Hodge structures. 
Before proving these results, let us calculate the 
cohomology of local systems on the total space of a
fibration with the fibers homotopy equivalent to $F=(\CC^*)^k$.
Let $\pi: T^* \rightarrow E$ be such a locally trivial fibration
over a manifold $E$ 
for which the associated $\CC^k$ bundle $T \rightarrow E$ is 
a direct sum of $k$ line bundles $L_1,...,L_k$. 
Denote by $c^1_1,...,c^k_1 \in H^2(E)$
the first Chern classes of $L_1,...,L_k$ respectively. 
Let $\rho_E$ be a unitary local system on $E$.
Consider the homomorphisms: 
\begin{equation}
\kappa^{a,b}: H^a(E,\L_{\rho_E}) \otimes \Lambda^b(H_1(\CC^*))
\rightarrow H^{a+2}(E,\L_{\rho_E}) \otimes \Lambda^{b-1}(H_1(\CC^*))
\end{equation}
 given by:
\begin{equation}
\kappa^{a,b}(\beta \otimes \alpha_1 \wedge ...\wedge \alpha_k)
=\sum_i \beta \cup c^j_1 \otimes \alpha_1 \wedge ...\hat \alpha_j ...
\wedge \alpha_k
\end{equation}
and let 
\begin{equation}
     K^{a,b}={\rm Ker} \kappa^{a,b}/{\rm Im } \kappa^{a-2,b+1} 
\end{equation}

Denote by  
$Im \pi_1(F)$ the image of the homomorphism $\pi_1(F) 
\rightarrow \pi_1(T^*)$ and 
let $\rho: \pi_1(T^*) \rightarrow \CC^*$ 
be a unitary representation of the fundamental group of $T^*$.
The exact sequence:
\begin{equation}
0 \rightarrow Im \pi_1(F) \rightarrow \pi_1(T^*) 
\rightarrow \pi_1(E) \rightarrow 0
\end{equation}
shows that a representation $\rho$ has trivial restriction on $Im \pi_1(F)$ 
if and only  if $\rho=\pi^*(\rho_E)$ for some unitary representation 
$\rho_E$ of $\pi_1(E)$.

\begin{lemma}\label{charbundle}
 For any non-negative $i$ one has the following:
\end{lemma}
\begin{equation}\label{charbundle1}
{\rm dim} H^i(T^*,\L_{\rho})=
\begin{cases}0 & \text{if \ \ $\rho \vert_{{\rm Im} \pi_1(F)} \ne 1$} 
\\ \sum_{a+b=i}{\rm dim} K^{a,b} & 
\text{if $\rho \vert_{Im \pi_1(F)}=1 \ {\rm or \ equivalently}\ 
\rho=\pi^*(\rho_E).$ }
\end{cases}
\end{equation}

\begin{proof} (of lemma \ref{charbundle}) 
(\ref{charbundle1}) is clearly the case when $E$ is a point 
which also
yields that the fiber of $R^q\pi_*(L_{\rho})$ for 
arbitrary $E$ is either $\Lambda^q(H^1(F,\CC))$ (if $\rho_E$ 
is trivial) or zero (otherwise).   
For a more general $E$, the lemma follows from the Leray spectral sequence:
$H^p(E,R^q\pi_*(\L_{\rho})) \Rightarrow H^{p+q}(T^*,\L_{\rho})$ which 
degenerates in the term $E_3$. Indeed if $\rho$ is not induced 
from $E$ then $E_2^{p,q}=0, p, q \ge 0$ by (\ref{charbundle1}). 
If $\rho=\pi^*(\rho_E)$
for a local system $\rho_E$ on $E$ then, since the local systems 
$R^q\pi_*(\CC)$ are trivial 
i.e. $E_2^{p,q}=H^p(E,\rho_E) \otimes \Lambda^q(H_1(F,\CC))$,  
the degeneration can be seen using 
the compatibility property of differentials of a spectral 
sequence with products and vanishing of $d_i, i>2$ on $E_*^{0.1}$.
One arrives at the formula (\ref{charbundle1}) for $E_3=E_{\infty}$ 
since $d_2^{a,b}=\kappa^{a,b}$ as one can see using multiplicativity 
and the well-known identification of $d_2$ in the Leray spectral sequence
of a circle bundle with the cup product with the first Chern class
of the associated line bundle.
\end{proof}

\begin{remark}\label{torusbundlecase} {\it Local systems on $\CC^*$-fibrations
over quasiprojective manifolds.}
Let $T^*$ be a 
$(\CC^*)^r$-fibration over 
$E$ such that the associated $\CC^r$-bundle $T$ is a direct sum of 
holomorphic line bundles and let:
\begin{equation}\label{translatedsbgstratum}
\S^n_l(T^*)=\{\chi \in \Char \pi_1(T^*) \vert H^n(T^*, \L_{\chi}) \ge l \}
\end{equation}
If the base $E$ is quasiprojective, then it follows from 
lemma \ref{charbundle} and \cite{Arapura} 
that $\S^n_l(T^*)$ is a finite union of 
translated subgroups. 
Indeed, $\S^n_l(T^*)$ is an invariant of the homotopy type and,  
since $T^*$ is homotopy equivalent to the complement in the total space of $T$
to union of total spaces of rank $r-1$ subbundle of $T$
which is also quasiprojective, one can apply \cite{Arapura}. 

Similarly, let $\rho$ be a unitary local system  
on a the total space of a fibration $T^*$. Assume again that
 the associated $\CC^r$ bundle $T$ is a
 direct sum of line bundles. Then the cohomology $H^p(T^*,\L_\rho)$ 
support the canonical $\CC$-mixed Hodge structure defined has follows. 
By lemma \ref{charbundle} we can assume that $\rho=\pi^*(\rho_E)$.
Consider the total space of bundle $T=\oplus_{i=1}^{i=r} L_i$.  The complement 
to the union of divisors $\oplus_{i \ne j} L_i \ \ (j=1,...,r)$ 
has the homotopy 
type of $T^*$ and the mixed Hodge structure on $H^p(T^*,\rho_E)$ 
is obtained from the mixed Hodge structure constructed in 
\cite{Timm}, \cite{Arapura} via the identification 
\begin{equation}H^p(T^*,\pi^*(\rho_E))=
H^p(T-\cup_j \oplus_{i \ne j} L_i,\pi^*( \rho_E))
\qedhere
\end{equation}
In the case when $\rho$ is trivial 
this is compatible with the standard mixed Hodge structure 
of a punctured neighbourhood. 
\end{remark}

\bigskip

The proofs of the theorems \ref{translatedsbg} and \ref{localmhs}
use the following result on degeneration of Mayer Vietoris 
spectral sequence generalizing the case of projective manifolds 
discussed in \cite{CG} to quasiprojective 
situation with local systems.

\begin{lemma}\label{degenerationMV}
Let $\bar X=\cup \bar X_i$ be a union of projective 
manifolds having only normal crossings and let 
${\bf D}_i=\cup_j D_{i,j} \subset \bar X_i$ be a normal crossing divisor 
such that $D_{i,j} \cup \bar X_i \cap_{k \ne i} X_k$ is 
a normal crossing divisor
 on $\bar X_i$ as well. Denote by $X_i=\bar X_i-\cup D_{i,j}$, 
$X_{i_0,..i_q}=X_{i_0} \cap ... \cap X_{i_q}$ and 
$X^{[q]}=\coprod X_{i_0,...,i_q}$.
Let $\rho: \pi_1(\cup X_i) \rightarrow \CC^*$ be a local 
system on $X$ and $\rho_{i_0,..,i_q}$ be the local system on 
$X_{i_0,...,i_q}$ induced 
by $\rho$. Let $H^p(X^{[q]},\L^{[q]})=\oplus_{i_0 <....<i_q} 
H^p(X_{i_0,...,i_q},\L_{i_0,...,i_q})$.
This data determines a spectral sequence:
\begin{equation}\label{firstterm}
E_1^{p,q}= H^p(X^{[q]},\L^{[q]}) \Rightarrow H^{p+q}(X,\L)
\end{equation}
degenerating in term $E_2$.
\end{lemma}

\begin{proof} Let $\M_{i_0,...,i_q}$ be the Deligne's extension
 (cf. \cite{delignediff}) 
corresponding to the local system $\L_{i_0,...,i_q}$ and 
$\M^{[q]}$ be the corresponding bundle on $X^{[q]}$. 
The spectral sequence in this lemma is the spectral 
sequence of the double complex
\begin{equation}\label{doublecomplex}
A^{\bullet,\bullet}=A^p(X^{[q]}(log),\M^{[q]},(\nabla, \delta))
\end{equation} 
with components being $C^{\infty}$ logaritmic forms 
on $\bar X^{[q]}(=\coprod \bar X_{i_0,...i_q})$
 with the poles along the intersections of components
of ${\bf D}_i$ with $\bar X_{i_0,...i_q}$ 
(twisted deRham complex for  
the Deligne's extension $\M^{[q]})$. The differentials respectively 
are the differental 
$\nabla: A^{p,q} \rightarrow A^{p+1,q}$ of the connection 
and the differential 
$\delta: A^{p,q} \rightarrow A^{p,q+1}$ 
which takes $\omega \in A^{p,q}$ having components $\omega(i_0,..,i_q) \in 
A^p(X_{i_0,...,i_q}(log),\L \vert_{X_{i_0,...,i_q}})$
to $\delta(\omega)(i_0,...i_{q+1})=
\sum (-1)^k\omega(i_0,...\hat i_k,...i_q) \vert_{X_{i_0,...,i_{q+1}}}$.

Clearly the term $E_1^{p,q}$ of the spectral sequence of 
the double complex 
(\ref{doublecomplex}) coincides with expression in (\ref{firstterm}).
The abutment of this spectral sequence is $H^*(X,\L)$ since the single 
complex $(A^{\bullet,\bullet},\nabla+\delta)$ is a locally free 
{\it acylic} resolution of locally constant vector bundle 
of the local system $\L$. Acyclicity 
follows from the local calculation which identifies for a sufficiently 
small open subset $U \subset X$ the group $H^n(U,\L)$ with the 
abutement of the spectral sequence having $E_2^{q,p}=H^q_{\delta}H^p_{\nabla}$
and vanishing of the latter: indeed for $U$ over which $\L$ is not trivial
one has $H^p_{\nabla}=0$ and for $U$ over which $\L$ is trivial 
$E_2^{0,q}$ is the cohomology of the simplex (cf. \cite{CG}).

To show the $E_2$-degeneration 
of (\ref{firstterm}) first let us consider the case when 
none of the eigenvalues of the holonomy of the connection about the components
of ${\bf D}_i \cap X^{[q]}$ is equal to one. In this case, for each  component
$X_{i_0,...,i_q}$ one has: $$H^i(X_{i_0,...,i_q},\L_{i_0,...,i_q})=
H^i(\bar X_{i_0,...,i_q},j_*\L_{i_0,...,i_q})=
H^i_{(2)}(X_{i_0,...,i_q},\L_{i_0,...,i_q})$$
 (cf. 
\cite{Timm},\cite{Timm2}; 
here $j$ is the 
embedding $X_{i_0,...,i_q} \rightarrow \bar X_{i_0,...,i_q}$ and 
$L^2$-cohomology are with respect to a complete metric 
on $X_{i_0,...,i_q}$ 
asymptotic to the Poincare metric).
Using the twisted analog of $\partial \bar \partial$ lemma 
(cf. \cite{Zucker})
as in 
compact case (\cite{CG}) one obtains degeneration 
in ths case when all monodormies of $\L_{i_0,...i_q}$ 
are non trivial.

General case now can be deduced by induction over the number 
of components $\bigcup_i{\bf D}_i$ for which $\L$ has non trivial monodormy.
Let $X'=(X_1-D) \cup 
\bigcup_{i \ge 2} X_i$ i.e. $X'$ is a quasiprojective normal crossing
having one more compactifying component than does $X$. 
Let $\L$ be a local system on $X'$ such that the holonomy about $D$ 
is trivial. The connection matrix of the flat connection 
corresponding to 
the local system $\L \vert_{X_{1,i_1,...,i_q}}$ has an integer 
residue along $D$ since the holonomy is trivial along $D$. 
Hence in the Deligne's 
extension this residue is equal to zero along $D$. 
In particular we have a well defined holomorphic 
connection on $X$ and its restriction to $D$.
Then one has 
\begin{equation}\label{sequenceofspectralsequences}
0 \rightarrow A^{p}(log)(X^{[q]}, \M)
\rightarrow A^{p}(log)(X'^{[q]}, \M)
 \rightarrow A^{p-1}(log)(D^{[q]}) \rightarrow 0
\end{equation}
where the last term is the log complex of $D^{[q]}=
D \cap X_{\ne 1}^{[q-1]}$ with $X_{\ne 1}^{[q-1]}$ 
is disjoint union of intersections of irreducible components of 
$X$ different from $X_1$ and the last map is the 
resudue map of log-forms along $D$.

Recall that the degeneration of a spectral sequence of  
a differential graded algebra $A^{\bullet}$ filtered 
by a decreasing filtration $F^p(A)$ 
in term $E_k$ is equivalent to condition:
\begin{equation}\label{degeneracy}
   F^p(A^{\bullet}) \cap dA \subset dF^{p-k+1}(A^{\bullet})
\end{equation}
(cf. \cite{morgan} Lemma (1.5) p.144 and \cite{delignehodge},  
(1.3.2) and (1.3.4)) 
One sees directly, that conditon (\ref{degeneracy}) 
for the endterms of (\ref{sequenceofspectralsequences})
yields it for the middle terms as well.

\end{proof}

\begin{proof} (of theorem \ref{translatedsbg})
{\it Step 1. Identification of $\X-\D$ with a union of torus bundles.} 
Let $\bigcup_{i \in I} E_i$ be the exceptional locus of a log resolution 
of the pair:
 $$(\tilde \X,\bigcup_{i \in I} E_i \cup \cup_{j=1}^{j=r}\tilde \D_j) 
\rightarrow (\X,\cup D_j)$$ where 
$\tilde D_j$ are the proper 
preimages of the components $D_j$. 
Recall that a resolution of a pair $(X,D)$ where $X$ is a normal 
variety is a morphism $f: Y \rightarrow X$ such that the union 
of the exceptional locus of $f$ and the proper preimage of $D$ in $Y$ 
is a normal crossing divisor. Given $(X,D)$ such a morphism $f$ always exist 
(cf. \cite{kollarmori}).
Let $\partial T(\cup E_i)$ be the boundary of a regular neighborhood
$T(\cup E_i)$ of $\cup E_i$ in $\tilde \X$. 
We have the identification (homotopy equivalence)
$$\partial T(\cup E_i) -\partial T(\cup E_i) \cap \cup \tilde D_j
=\X-\cup D_j$$
The (open) manifold 
$\partial T(\cup E_i) -\partial T(\cup E_i) \cap \cup \tilde D_j$
can be constructed inductively as a union of tori bundles 
over quasiprojective varieties 
intersecting
along unions of tori bundles lower dimension. More precisely, 
consider the stratification of $\bigcup E_i$ in which 
each stratum is a connected component of a set 
consisting of points belonging to exactly $l$ components of this 
union. Complement to $\cup D_j$ in an intersection 
of the boundaries $\partial T(E_i)$
of small regular neighborhoods in $\tilde \X$ of 
$E_i$ containing this stratum 
is  a fibration 
with the fiber being the torus $(S^1)^l$.
Here $l$ is the number of components $E_i, i \in I$ containing this stratum.
Each such torus fibrations has a compact base if and only if the 
corresponding component $E_i$ of the exceptional locus has empty intersection 
with the proper preimage of $D \subset \X$. One has a locally trivial 
torus fibrations over the complement in an intersection of  components $E$ 
having non empty intersection with the proper preimage of $D$ due to 
normal crossing condition on union of the exceptional locus and 
the proper preimage of $E$ (note that we do not need assumption 
that $\X$ has an isolated singularity).

{\it Step 2. Translated subgroup property for each stratum.}
For each stratum $S$, the total space $T^*(S)$ of this $(S^1)^l$-fibration 
is a subset of 
$\partial T(\cup E_i) -\partial T(\cup E_i) \cap \cup \tilde D_j$. 
The collection of characters in 
${\bf T}={\rm Char} (\pi_1(\partial T(\cup E_i) 
-\partial T(\cup E_i) \cap \cup \tilde D_j))$,
which when restricted on $T^*(S)$ 
yields a character with the corresponding local system having 
$H^n(T^*(S),\L_{\chi}) \ge l$, is a union of translated subgroups 
in $\bf T$ (cf. (\ref{translatedsbgstratum}) or, since 
$S$ is quasiprojective, apply the remark after 
proof of lemma \ref{charbundle}).

{\it Step 3. Degeneration of Mayer-Vietoris spectral sequence 
for torus bundles.}
The cohomology $H^n(\partial T(\cup E_i) 
-\partial T(\cup E_i) \cap \cup \tilde D_j,\L_{\chi})$
is the abutment of a Mayer Vietoris 
spectral sequence:

\begin{equation}\label{mayer}
   E_1^{p,q}=H^p(A^{[q]}, \L_{\chi \vert_{A^{[q]}}}) 
\Rightarrow H^{p+q}(\bigcup A_i \L_{\chi})
\end{equation}
where $A^{[q]}$ is a torus bundle over a stratum of the above stratification.
This 
spectral sequence has trivial differentials $d_i$ for $i \ge 2$.
Indeed, the degeneration of the spectral sequence (\ref{mayer}) 
in the case when $A^{[0]}$ are quasiprojective was shown in lemma 
\ref{degenerationMV}.
The case when $A^{[q]}$ are tori bundles over 
quasiprojective manifolds follows from formulas in lemma \ref{charbundle}.

{\it Step 4. End of the proof.} 
Now let us consider the collection of quasi-affine subsets   
$\cal T$ of the torus  $\bf T$ with 
each quasi-affine 
subset being a finite union of subgroups of $\bf T$ with a removed 
collection of translated subgroup of $\bf T$ (possibly empty).
Each collection of local systems with fixed dimension $H^n(A, L_{\chi})$,
in the case 
when $A$ is quasiprojective, belongs to $\cal T$ and 
hence the collection $\chi$ for which each $E_2^{p,q}$ 
term in the spectral sequence (\ref{mayer}) has a fixed dimension 
also belongs to $\cal T$ as follows from Corollary 
\ref{restriction} (since $d_1$ in (\ref{mayer}) 
is the restriction map). This proves the theorem.
\end{proof}

\begin{proof} (of theorem \ref{localmhs} ) The main point
 is that, while making the calculation of 
the cohomology of local system  on resolution of $\X,\D$ 
using the isomorphism 
$H^i(\X-\D,\L_{\chi})=H^i(\tilde \X-(\bigcup_{i \in I} E_i \cup \tilde \D))$
where $\tilde \X,\cup_{i \in I} E_i \cup \tilde D$ is a resolution of  
pair $\X,\D$, one can replace $\bigcup_{i \in I} E_i \cup \tilde \D$
by the union of components $E_i$ such that the restriction of $\chi$
of $\pi_1(\partial T(E^{\circ}_i)-\tilde \D)$ is a pull back of a character 
of $\pi_1(E_i^{\circ}-\tilde \D)$ with respect to 
projection 
$\partial T(E^{\circ}_i)-\tilde \D \rightarrow E_i^{\circ}-\tilde \D$
(here $E_i^{\circ}$ is the locus of points of $E_i$ which 
are nonsingular points of $\bigcup E_i$;
the set of indices of $E_i$'s labeling these components will be 
denoted $I_{\chi}$).
In quasiprojective case this corresponds to isolating components with 
trivial holonomy of the connection. 
(cf. \cite{Timm},\cite{Arapura}). Moreover, in the case of ${\CC^*}^r$-bundles
we have the mixed Hodge structure
 as described in the remark \ref{torusbundlecase}.
Then all steps used in \cite{Timm} to derive the 
extension of \cite{delignehodge} with modifications as in local case of 
\cite{DH}, i.e. constructing the Mixed Hodge Complex yielding the 
Mayer Vietoris spectral sequence go through in our local case as well.

More precisely, let $E_{i_1,...,i_k}^{\circ}$ be the locus 
of points of $\bigcup_{i \in I} E_i$ which belong to the components 
$E_{i_1},...,E_{i_k}$ but do not belong to 
any other components of $E$ (recall that 
$I$ denotes the set of indices of components
of exceptional set). Let 
$\partial T(E_{i_1,...,i_k}^{\circ})-\tilde \D$ 
be  $T(E_{i_1} \cap ....\cap E_{i_k})-\bigcup E_i$
where $T(E_{i_1} \cap ....\cap E_{i_k})$ is a 
tubular neighborhood in $\X$. One has 
\begin{equation}
    \partial T(\bigcup E_i)-\D=\bigcup_{(i_1,...i_k) \subset I} 
\partial T(E_{i_1,...,i_k}^{\circ})-\tilde \D
\end{equation}
For each character $\chi \in \Char \pi_1(\X-\D,\L)$, let 
$\chi_{\partial T(E_{i_1,...,i_k}^{\circ})-\tilde \D}$ be the 
induced character of $\pi_1(\partial T(E_{i_1,...,i_k}^{\circ})-\tilde \D)$.
Let, as above, $I_{\chi}$ be the collection of components of exceptional 
divisor for which $\chi_{\partial T(E_i)}$ is a pull back 
of a character $\chi_{E_i^{\circ}-\D \cap E_i^{\circ}}$
(i.e., which restriction of the boundary of transversal to 
$E_i$ circle is the  trivial character: cf. (\ref{charbundle});
we shall call such components {\it fiber $\chi$-trivial}).
The characters of $\pi_1(E_{i_1,...,i_k}^{\circ}-E_{i_1,...,i_k}^{\circ}
\cap \D)$, i.e.,
 classes in $H^1(E_{i_1,...,i_k}^{\circ}-E_{i_1,...,i_k}^{\circ}
\cap \D,\CC^*)$ obtained
from the characters of $\tilde \X-\bigcup_{i \in I} E_i \cup \tilde \D$
(and hence compatible with restrictions),
 define as a result of degeneration of Mayer Vietoris spectral sequence,
the characters of $\pi_1(\bigcup_{i \in I_{\chi}} E_i)$ 
and hence the local system on $T(\bigcup_{i \in I_{\chi}} E_i)$ 
which we also denote $\L_{\chi}$.

Two Mayer Vietoris spectral sequences:

\begin{equation}\label{mayer1}
  E_1^{p,q}=H^p(\bigcup_{(i_1,...,i_{q-1}) \in I} 
\partial T(E_{i_1,...,i_{q-1}}^{\circ})-\tilde \D, 
\L_{\chi_{T(E_{i_1,...,i_{q-1}}^{\circ})-\tilde \D}})
\end{equation}
$$\Rightarrow H^{p+q}(\partial T(\bigcup E_i)-
\partial T(\bigcup E_i) \cap \D,\L_{\chi})
$$

and 

\begin{equation}\label{mayer2}
  E_1^{p,q}=H^p(\bigcup_{(i_1,...,i_{q-1}) \in I_{\chi}} 
\partial T(E_{i_1,...,i_{q-1}}^{\circ})-\tilde \D, 
\L_{\chi_{T(E_{i_1,...,i_{q-1}}^{\circ})-\tilde \D}})
\end{equation}
$$\Rightarrow H^{p+q}(\partial T(\bigcup_{i \in I_{\chi} E_i})-
\partial T(\bigcup E_i) \cap \D,
\L_{\chi_{\partial T(\bigcup_{i \in I_{\chi} E_i})-
\partial T(\bigcup E_{i \in I_{\chi} E_i}) \cap \D}}))
$$
yield that it is enough to show the existence of the 
Mixed Hodge structure on the cohomology 
of local systems of the boundary of the neighborhood only 
for those components 
for which the restriction of $\chi$ on the boundary 
is the pullback from this component. 

In order to describe the cohomological $\CC$-Hodge complex 
$(A,F,\bar F,W)$ on $\X-\D$ with cohomology being the 
cohomology of the abutment of the spectral sequence (\ref{mayer2})
first let us describe the complex for each intersection 
$\partial T(E_{i_1,...,i_{q-1}}^{\circ})-\tilde \D$ and the 
local system $\L_{\pi^*(\chi)}$ pulled back to it from the base. 
If $\TT(E_{i_1,...,i_{q-1}})$ is defined as the projectivization 
of the total space
of the direct sum of the normal bundle to $E_{i_1,...,i_{q-1}}$ in 
$\X$ with $\O_{E_{i_1,...,i_{q-1}}}$
and $\iota_{i_1,...,i_{q-1}}$ is the embedding of the neighbourhood 
of $E_{i_1,...,i_{q-1}}$ in $\TT(E_{i_1,...,i_{q-1}})$
then 
$\TT(E_{i_1,...,i_{q-1}})$ is a compactification of the total 
space of the normal bundle of $E_{i_1,...,i_{q-1}}$ in $\X$. 
$\TT(E_{i_1,...,i_{q-1}})$ is also a compactificaiton of 
$\CC^*$-fibration over $E_{i_1,...,i_{q-1}}^{\circ}-
E_{i_1,...,i_{q-1}}^{\circ} \cap \tilde \D$ 
and the complement in it to the total space of this fibration 
is a normal crossing divisor which we denote as 
$D_{i_1,..,i_{q-1}} \subset \TT(E_{i_1,...,i_{q-1}})$.  
This fibration is 
homotopy equivalent to the punctured neighbourhood of 
 $E_{i_1,...,i_{q-1}}^{\circ}-\tilde \D$ in $\X$.
We can use the just described 
compactification to construct (following \cite{Arapura}) 
log-complex 
associated with the  (Deligne's extension) $\bar V$ of the connection
with the holonomy given by $\rho$. We consder trifiltered real 
analytic log-complex:
\begin{equation}\label{onecomponentcomplex}
A^{\bullet}_{i_1,...,i_q}=
(A^{\bullet}_{\TT(E_{i_1,...,i_{q-1}})}(\log D_{i_1,..,i_{q-1}})
\otimes \bar V,F,\bar F,W)
\end{equation}
used in \cite{Arapura} and generalizing the complex considered
in \cite{Navarro} in the case $\bar V$ is trivial and providing
supporting conjugation analog on the 
log-complex with connection constructed in \cite{Timm}.
The log-complex (\ref{onecomponentcomplex}) is quiasiisomorphic to 
direct image of its restriction on punctured neighbourhood 
$T(E_{i_1,...,i_{q-1}})$ (i.e. with deleted zero-section) 
which we denote 
\begin{equation}
(\tilde A^{\bullet}_{i_1,...,i_q} \otimes \bar V,F,\bar F,W) 
\end{equation}
Here $F,\bar F,W$ are the filtrations induced by corresponding 
filtrations of (\ref{onecomponentcomplex})
One has the maps $\delta:\tilde A^{\bullet}_{i_1,...,i_q} \rightarrow 
\tilde A^{\bullet}_{i_1,...,i_q,i_{q+1}}$ induced by inclusion 
of punctured neighbourhoods which allow to form differential graded complex 

\begin{equation}
\tilde A^{\bullet,\bullet}=\oplus_{i_1,..,i_{q-1}} \tilde 
A^{\bullet}_{i_1,....,i_{q-1}}
\end{equation}

The total complex $T^n=\oplus_{i+j=n} \tilde A^{i,j}$ 
with the usual differential $d+(-1)^j\delta_j$  and 
filtrations $F^pT^n=\oplus_{i \ge p,i+j=n} \tilde A^{i,j}, 
\bar F^p=\oplus_{j \ge p
,i+j=n}, W_kT=\oplus W_{k+i} \tilde A^{i,j}$ 
(cf. \cite{Navarro},\cite{Arapura1}) is a $\CC$-mixed Hodge complex
calculating the abuttment of the spectral sequence (\ref{mayer2}) and hence
(\ref{mayer1}) i.e. yielding the Mixed Hodge structure on the 
cohomology of the abutment. Compatibility with holomorphic maps 
follows form the corresponding compatibility in quasiprojective case.

\end{proof}

\begin{remark}
Assume that one has an embedding $(\X,\D) \subseteq (\bar \X,\bar \D)$ 
such that $(\bar \X,\bar \D)$ is quasiprojective
and such that $\chi$ is a pullback of a character of 
$\pi_1(\bar \X-\bar \D)$.
For a germ on an isolated 
non-normal crossing in $\CC^n$, the extension of each irreducible 
component  $\D$ to an irreducible hypersurface in $\CC^n$ yields
such $(\bar \X, \bar \D)$.
The mixed Hodge structure of the theorem \ref{localmhs} on $H^i(\L_{\chi})$
where $\L_{\chi}$ is the local system on the boundary 
of a punctured neighbourhood of union of fiber $\chi$-trivial 
exceptional divisors minus $\D$ can be constructed using methods of \cite{DH}
(by above, these cohomology coincide with 
$H^i(\X-\D,\L_{\chi})$).
Indeed, such a punctured neighbourhood is the intersection 
of regular neighbourhood of union of fiber $\chi$-trivial divisors
and $\bar X-\bar D$. Both spaces support 
the canonical ($\CC$-) mixed Hodge structure.  
In particular,  the cohomology of such a local system on the boundary of
a neighborhood of union of fiber $\chi$-trivial exceptional divisors
can be 
calculated using hypercohomology of the 
mapping cone of the following bifiltered complexes of sheaves
corresponding to presentation 
$T(\bigcup_{i \in I_{\chi}} E_i-\D \cap \bigcup_{i \in I_{\chi}} E_i)$
as an intersection of 
a quasiprojective variety containing $\tilde \X-\tilde \D$
(i.e., the resolution of singularities $\tilde {\bar \X},
 \tilde {\bar \D}$
of $\bar \X,\bar \D$
in which preimage of  $\X,\D$ is $\tilde \X, \tilde D$)
and a regular neighborhood of $\bigcup_{i \in I_{\chi}} E_i-\D 
\cap \bigcup_{i \in I_{\chi}} E_i$
(this construction in a standard way can be upgraded
to the level of mixed Hodge complexes).
The first complex has as its hypercohomology the the cohomology of the 
regular neighborhood of $T(\bigcup_{i \in I_{\chi}} E_i)-\D$ and is given 
by:

\begin{equation}
\Omega^{\cdot }_{\bigcup_{i \in I_{\chi}}E_i}(log \D)
=\oplus_{p+q=n} 
(\pi_{q})_*\Omega^{p}_{E^{[q]},i \in I_{\chi})}({\rm log} (\D \cap
 \cup E))
\otimes \V_{\chi}
\end{equation}
i.e., the single complex associated to the double
complex with components being the push forward to 
$\bigcup_{i \in I_{\chi}} E_i$ of log forms with the 
values in the induced (from $\tilde {\bar \X}$)
the Deligne's extension on all  
$q+1$ fold intersections $E^{[q]}$
of $E_i, i \in I_{\chi}$ relative to divisors on 
connected components of $E^{[q]}$ induced by intersections with
$\D$ and with differential given by the log-connection.

The second complex, 
having as its hypercohomology the cohomology of the complement 
to $\bigcup_{i \in I_{\chi}} E_i \cup \D$ in the resolution 
of singularities of quasiprojective variety $\bar X$,
is given by:

\begin{equation} \Omega^{\cdot}_{\tilde{\bar \X}}({\rm log} 
\cup_{i \in I_{\chi}} E_i \cup \D)
\otimes \V_{\chi}
\end{equation}
The Hodge filtration and weight filtrations
are given in the usual way (i.e., by truncation and 
 by considering forms with vanishing $m$-residues  
respectively).
\end{remark}

\noindent The construction in the theorem yields:

\begin{corollary}\label{functoriality}
 The above mixed Hodge structure is functorial in the 
following sense. Consider the homomorphism:
$$h^k(f): H^k(\widetilde {{\cal X}-{\cal D}},{\cal V}) \rightarrow 
H^k({{\cal X}-{\cal D}},f_*({\cal V}))$$
induced by an unbranched covering 
$f: \widetilde {{\cal X}-{\cal D}} \rightarrow {{\cal X}-{\cal D}}$
(here $f_*({\cal V})$ is the direct image of the local system, 
corresponding to the induced character of subgroup 
$\pi_1(\widetilde {{\cal X}-{\cal D}})$ of 
$\pi_1({{\cal X}-{\cal D}})$). Then $h^k(f)$ is a morphism of 
mixed Hodge structures.
\end{corollary}

The above proof of theorem \ref{translatedsbg} allows one to deduce
the following property of the mixed Hodge structures
 described in theorem \ref{localmhs}:

\begin{corollary} Let $\X,\D$ be as above. Then the subset
of a fundamental domain in the universal cover $\U$ 
of $\Char^u(\pi_1(\X-\D)$ given by:
\begin{equation}
\{u \in \U \vert \dim Gr_F^pH^{n}(\X-\D, V_{{exp(u)}}) \ge l\}
\end{equation}
is a finite union of polytopes.
\end{corollary}

\begin{proof} It follows from the degeneration 
of the Mayer-Vietoris spectral sequence mentioned in the proof 
of theorem \ref{translatedsbg}, the compatibility of $d_1$ in it with 
the mixed Hodge structures and the theorem in the previous section.  
\end{proof}

This extends the results of \cite{alexhodge} 
on Hodge decomposition of characteristic varieties 
of germs of plane curves where also 
examples of such ``Hodge decomposition'' of $S^{1,p}_l$ in this 
case are given.

\begin{remark} The construction of the mixed Hodge structure
 in the local case also can be extended
to the case of cohomology of a unitary local system, i.e. to the 
context of \cite{Timm}. The theorem \ref{translatedsbg} 
can be modified in an obvious way 
to include a twisting by a higher rank local system as in theorem
 \ref{quasiproj}.
\end{remark}

\section{Isolated non-normal crossings}\label{inncsection}

Now we shall apply the results of previous 
section to the case of isolated non-normal crossings.
First recall the following, already mentioned in the introduction:

\begin{dfn}\label{localinnc}
(cf. \cite{innc}, \cite{DL})
 An isolated non-normal crossing (INNC) is a pair $(\X,\D)$ 
where $\X$ is a germ of a complex space $\X$ having $\dim \X-2$-connected link
 $\partial B \cap \X$ where $B$ is a small ball about $P \in \D$
\footnote{for example a germ 
of a complete intersection with isolated singularity.} 
and where  
$\D$ is a divisor on $\X$ which has only normal 
crossings at any point of $\X-P$. 
\end{dfn}

The global isolated-non normal crossings divisors 
were considered in \cite{ample}
where the homotopy groups of the complement were related to 
the local invariants which are certain polytopes. The goal 
of this section is to relate them 
to the polytopes discussed in this paper.

For the local INNCs as in \ref{localinnc},
one has the following homotopy vanishing theorem 
(cf. \cite{innc}. Th.2.2
and \cite{DL}):

\begin{theorem}\label{vanishing}
 If $(\X,\D)$ is an INNC, $\dim \X=n+1$ and $r$ is the 
number of irreducible components in $\D$ then:
$$\pi_1(\X-\D)=\ZZ^r, \ \ \pi_i(\X-\D,x)=0 \ \ {\rm for} \ 2 \le  i \le n-1$$
\end{theorem}

The main invariant of local INNCs is $\pi_n(\X-\D)$ 
considered as a $\pi_1(\X-\D)$-module. It follows from the theorem 
\ref{vanishing} that this homotopy group is isomorphic 
as a $\ZZ[\pi_1(\X-\D)]$-module to the 
homology of the infinite abelian cover $H_n(\widetilde {\X-\D},\ZZ)$.
The $\pi_1(\X-\D)$-module structure on the latter is given by 
the action of the fundamental group on the homology of 
the universal cover.

\begin{dfn}(cf. \cite{abcov}, \cite{topappl})
 $l$-th characteristic variety of $(\X,\D)$
 is the reduced support of 
$\pi_1(\X-\D)$ module $\Lambda^l (\pi_n(\X-\D) \otimes_{\bf Z} \CC)$. 
In other words:
\begin{equation}
S_l(\X,\D)=\{ \wp \in \Spec \CC[\pi_1(\X-\D)]={\CC^*}^r \ \ 
\vert \ \ 
(\Lambda^l (\pi_n(\X-\D) \otimes_{\bf Z} \CC))_{\wp} \ne 0
\}
\end{equation}
(here subscript $\wp $ denotes the localization in the prime ideal $\wp$).
\end{dfn}

INNCs are reducible analogs of isolated singularities corresponding 
to the case $r=1$. The 
homotopy vanishing in theorem \ref{vanishing}
is equivalent to the 
Milnor's theorem asserting that the Milnor fiber of an isolated 
singularity of a hypersurface having dimension $n$ is $(n-1)$ connected
 (cf. \cite{milnor}).
Indeed, by Milnor's fibration theorem one has a locally trivial 
fibration: $\phi: \X-\D \rightarrow S^1$. Clearly, the identities of 
the theorem \ref{vanishing} are equivalent to the requirement that 
the fiber $F$ of $\phi$  
satisfies: $\pi_i(F)=0\ \ 0 \le i \le n-1$. Moreover, in this case (i.e. when 
$r=1$)
the variety $S_1 \subset \CC^*$ is the collection of eigenvalues of the 
monodromy of Milnor fibration $\phi$. 

Note that since the universal cover $\widetilde {\X-\D}$ of $\X-\D$  
is also a covering space of the Milnor fiber and since the latter 
has the homotopy type of an $n$-complex, it follows that   
the universal cover of $\X-\D$ has the homotopy type of a bouquet of 
$n$-spheres (cf. remark 4.7 in \cite{innc}).
  
The characteristic varieties are equivalent to the loci considered 
in the theorem \ref{translatedsbg}:

\begin{prop}(\cite{DL},\cite{innc})
A local system $\chi \in \Char (\pi_1(\X-\D)), \chi \ne {\bf 
1}$ corresponds to a point $S_l \subset \Spec\CC[\pi_1(\X-\D)]$ after the
 canonical identification $\Char (\pi_1(\X-\D))=\Spec\CC[\pi_1(\X-\D)]$
if and only if 
  $$\dim H^n({\cal L}_{\chi}) \ge k$$
i.e. the jumping loci of the cohomology of local systems coincide with 
the supports of the homotopy group $\pi_n$.
\end{prop}

As result we obtain the following:

\begin{corollary}
The components of characteristic variety are translated 
subgroups by torsion points.  
\end{corollary}

This corollary can be viewed as a generalization of the classical monodromy 
theorem since in the case $r=1$ this is equivalent to the claim that 
an eigenvalue of the monodromy operator is a root of 
unity.

Next we shall interpret the results of previous sections 
in terms of mixed Hodge theory of abelian covers.

As above, for 
a group $G$ acting on a vector space $V$ and a character $\chi \in \Char G$
we shall denote by $V_{\chi}$ the subspace $\{v \in V \vert g \cdot v
=\chi(g) v\}$. The following is a direct generalization 
of Prop. 4.5 and 4.6 in \cite{innc} on homology of branched and unbranched 
covers and the relation with cohomology of local systems in the case of 
curves discussed in \cite{abcov} \cite{eko}.

\begin{prop}\label{homologycovers}
 a) Let $X=(\X-\D) \cap \partial B$ (where $\partial B$ is 
a small ball about the point of non normal crossing of $\D$). 
Let $U_{\bf m} \rightarrow X \ \ ({\bf m}=(m_1,...,m_r))$ be the abelian 
covering corresponding to the homomorphism $\pi_1(X)=\ZZ^r \rightarrow 
\oplus \ZZ/m_i\ZZ$. For ${\bar \omega}=(...,\omega_i,...) \in {\CC^*}^r$ let 
\begin{equation}
f(\bar \omega,X)={\rm max}\{l \vert \bar \omega \in S_l( X) \}
\end{equation}
Then 
\begin{equation}
 H_k(U_{\bf m},\ZZ)=\Lambda^k(\ZZ^r) \ \ \ (1 \le k <n)
\end{equation}
\begin{equation}
{\rm rk} H_n(U_{\bf m},\QQ)=\sum_{\bar \omega, \omega_i^{m_i}=1}
f(\bar \omega,X)
\end{equation}

b) Let $V_{\bf m}$ be the branched covering space of $\X \cap \partial B$
branched over $\D \cap \partial B$. For a character 
$\chi \in \Char H_1(\X-\D,\ZZ)$ let 
\begin{equation}
I_{\chi}=\{ i \vert 1 \le i \le r \ {\rm and} \ \chi(\gamma_i^{j_i}) \ne 1
\ {\rm for \ some}\  j_i, 1 \le j_i <m_i \}
\end{equation}
where $\gamma_i \in H_1(\X-\D,\ZZ)$ is represented by the boundary of 
a small disk in $\X$ transversal to the component $D_i$ of $\D$.
Then $\chi$ can be considered as the character of 
$H_1(\X-\cup_{i \in I_{\chi}} D_i, \ZZ)$ and 
$$H_i(V_{\bf m},\ZZ)=0 \ \ \ 1 \le i \le n-1 $$
\begin{equation}
H_n(V_{\bf m},\ZZ)=\sum_{\chi \in \Char(\oplus_i {\ZZ/m_i\ZZ})}
f(\chi,\X-\cup_{i \in I_{\chi}} D_i)
\end{equation}
\end{prop}

If the components $D_i$ of $\D$ are zeros of holomorphic functions:
$f_i$ for $i=1,...,r$ respectively and $\X \subset \CC^N$ as above, then 
$V_{\bf m}$ has the realization as the link of complex space 
$\X_{\bf m} \subset \X \times \CC^r \subset 
\CC^{N}\times \CC^r= \{ ({\bf x},z_1,...,z_r) \vert
{\bf x} \in \CC^{N}, (z_1,...,z_r) \in \CC^r \}$ such that:

\begin{equation}
  z_1^{m_1}=f_1({\bf x}), . \ .\  .\  , z_r^{m_r}=f_r({\bf x}),  
\ \ \  {\bf x} \in \X
\end{equation} 
 
\begin{prop}\label{coverslocalsystems} 
Let $G=\oplus \ZZ/m_i\ZZ$ and let 
$\chi \in \Char G \subset \Char \ZZ^r$. The group $G$ acts on 
$H_n(U_{\bf m},\QQ)$ and 
\begin{equation}
H_n(U_{\bf m},\QQ)_{\chi}=H_n(\L_{\chi})
\end{equation}
If $\chi(\gamma_i^{j_i}) \ne 1$ for $1 \le i \le r$
where $\gamma_i,j_i$ as in the 
proposition \ref{homologycovers} (b) then there is canonical 
isomorphism:
\begin{equation}
H_n(V_{\bf m},\QQ)_{\chi}=H_n(\L_{\chi})
\end{equation}
\end{prop}

Next we have the following corollary of the functoriality 
of the mixed Hodge structure:

\begin{prop} The Mixed Hodge structure
 on the cohomology of an abelian cover 
with a finite Galois group $G$ 
of a complement ${\cal X}-{\cal D}$ to an INNC
determines the mixed Hodge structure on the cohomology of 
essential (cf. \cite{innc})
local systems of 
finite order via:
  $$h^{p,q}(\tilde {{\cal X}}_G)_{\chi}=
h^{p,q}_{\chi}(({\cal X}-{\cal D})_{G})=h^{p,q}({\cal L}_{\chi})$$
Here $({\cal X}-{\cal D})_{G}$ is the Galois cover with group $G$
and $\tilde {{\cal X}}_G$ is the corresponding branched cover.
\end{prop}

\begin{proof} 
We have the isomorphism 
$H^n(({\cal X}-{\cal D})_{G})_{\chi}=H^n({\cal L}_{\chi})$
(for example one can use a chain description of the cohomology of local 
systems cf.\cite{abcov} or the 
next section for twisted version 
of it) which by corollary \ref{functoriality} is compatible with both 
filtrations. The first equality follows from the assumption that 
$\chi$ is essential. 
\end{proof}

This yields the following calculation of the Hodge numbers 
of abelian covers:  

\begin{prop} 
$\dim Gr_F^pH^n(U_{\bf m})_{\chi}= {\rm max} \{ \ l \  \vert \chi=exp(u) 
\ where \ u \in S^{n,p}_l \}$. 
If $\chi$ is as in Proposition \ref{coverslocalsystems} then 
this is also equal to 
$\dim Gr_F^pH^n(V_{\bf m})_{\chi}$.
\end{prop}

In particular the function $f(m_1,...,m_r)=\dim Gr^p_FH^n(U_{\bf m})$
is polynomial periodic with the degrees of polynomial depending 
on the dimensions of the polytopes $S^{n,p}_l$. This result 
for the Hodge groups $H^{p,0}$ was obtained in \cite{budur}.

\section{Twisted Characteristic Varieties}

In this section we describe a 
multivariable generalization (cf. definition \ref{twistedcharvar}) 
of the twisted 
Alexander polynomials considered in \cite{CF} and relate it 
to the
subvarieties $S_{\rho}^{n,l}$ described in theorem \ref{quasiproj}.
It is similar to the way the characteristic varieties
generalize the Alexander polynomials (cf.\cite{abcov}).
The theorem \ref{roots} extends the 
cyclotomic property of the roots of Alexnader polynomials of 
plane algebraic curves (cf. \cite{duke}) to twisted case. 
We use theorem \ref{roots} 
to describe a group theoretical property of the fundamental groups 
which may not be shared by other groups. We also use twisted characteristic
variety to obtain information on the homology of non unitary local systems
(cf. Prop. \ref{proptwistedcharvar}).

Let $X$ be a finite CW complex and $\tilde X$ (resp. $\tilde X_{ab}$) 
denotes the universal cover of $X$, 
(resp. its universal abelian cover). As usual
 $\pi'_1(X)=\pi_1(\tilde X_{ab})$ 
be the commutator of the fundamental group of $X$. Let 
$C_*(\tilde X)$ be a chain complex of $\tilde X$ with 
$\CC$-coefficients. We shall assume
that each graded component is {\it a right} free $\CC[\pi_1(X))]$-module.
For a $\CC[\pi_1(X))]$-module $M$ we shall denote by 
$M_{ab}$ the $\CC[\pi_1'(X)]$-module 
obtained by restricting the ring $\CC[\pi_1(X)]$ to
$\CC[\pi_1'(X)]$. A unitary representation 
 $\rho: \pi_1(X) \rightarrow U(V)$ yields 
the left $\CC[\pi_1'(X)]$-module $V_{ab}$
(here $U(V)$ 
is the unitary group of a hermitian vector space $V$).
 The modules 
$C_*(\tilde X) \otimes_{\pi_1'(X)}  V_{ab}$ with the action 
of $g \in \pi_1(X,\ZZ)$ given by 
$g(c \otimes v)=
cg^{-1} \otimes g \cdot v$ form a complex of $\CC[H_1(X,\ZZ)]$
modules which results in the structure of a $\CC[H_1(X,\ZZ)]$-module 
on the homology of the abelian cover $H_i(\tilde X_{ab},V_{ab})$
with coefficients in the local system $V_{ab}$ 
(obtained by restriction of $\rho$ 
on the commutator of the fundamental group). 
Recall that we assume that $H_1(X,\ZZ)$ is torsion free
and hence $\Spec\CC[H_1(X,\ZZ)]$ is connected: the general case
requires only obvious modifications.

\begin{dfn}\label{twistedcharvar} The twisted by $\rho$ 
$l$-th characteristic variety (in degree $n$) (denoted as $\Ch_n^l$)
is the support of the module
$$\Lambda^l H_n(\tilde X_{ab},V_{ab})$$
i.e., the subset of $\Spec\CC[H_1(X,\ZZ)]$ consisting of  
prime ideals $\wp$ in $\CC[H_1(X,\ZZ)]$ such that localization 
of  $\Lambda^lH_n(\tilde X_{ab},V)$ at $\wp$ is non zero.
\end{dfn}

Let us in addition fix  a surjection 
$\epsilon: \pi_1(X) \rightarrow \ZZ$, and denote $\bar \pi=\Ker \ \epsilon$.
Let $\rho: \pi_1(X) 
\rightarrow U(V)$ be a unitary representation.
With such data one associates a twisted Alexander polynomial and
the homology torsion which both were   
studied in \cite{KL} \cite{CF} (cf. references there for the prior work).
The main geometric situations in  which these invariants  
appear are the case when $X$ is a
complement to a knot (cf. \cite{KL}) and the case of the complement to 
a plane algebraic curve in $\PP^2$ (cf. \cite{CF}).

Recall that the twisted Alexander polynomials and the torsion
corresponding to the 
data $(X,\epsilon,\rho)$ 
are defined as follows. Condisder the 
covering space $X_{\bar \pi}$ corresponding 
to the subgroup $\bar \pi$ and the local system $V_{\bar \pi}$ 
on the latter obtained by restricting $\rho$ to $\bar \pi$. 
The $k$-th Alexander polynomial $\Delta_k(t)$ is 
the order of the $\CC[\ZZ]$- module 
${\rm Tor}H_k(\tilde X_{\bar \pi},V_{\bar \pi})$ with the convention that 
the order of a trivial module is $1$ (cf. \cite{KL},  Def.2.3). 
The homology torsion is defined as:
\begin{equation}
\Delta_{\epsilon,\rho}(t)=\Pi_k \Delta_k(t)^{(-1)^{k+1}}
\end{equation}
In other words, $\Delta_k=\Pi a_{k,i}(t)$ where
${\rm Tor}H_k(\tilde X_{\bar \pi},V_{\bar \pi})
=\oplus_i (\CC[\ZZ]/(a_{k,i}(t)))$
($a_{k,i}(t)$ divides $a_{k,i+1}(t)$) is the cyclic decomposition.

If $X$ has a homotopy type of a 2-dimensional complex 
such that $H_1(X,\ZZ)=\ZZ$
(which is the case for the complements to irreducible 
algebraic curves in $\PP^2$ or knots in $S^3$)
and $\epsilon$ is the abelianization:
$\pi_1(X) \rightarrow H_1(X,\ZZ)$ (i.e $\bar \pi$ is the commutator 
$\pi_1'$), 
then the homology torsion is just ${{\Delta_1(t)} \over {\Delta_0(t)}}$ 
since the $\CC[\ZZ]$- module $H_k(\tilde X_{\pi_1'},V_{\pi_1'})$
is free for $k=2$ (and trivial for $k>2$). To be more specific, 
the chain complex
calculating $H_k(\tilde X_{\pi_1'},V_{\pi_1'})$ is given by:

\begin{equation}\label{milnortypeseq}
 0 \rightarrow C_2(\tilde X) \otimes_{\pi_1'} V_{\pi'}
\rightarrow C_1(\tilde X) \otimes_{\pi_1'} V_{\pi'} 
\rightarrow C_0( \tilde X) \otimes_{\pi_1'} V_{\pi'} 
\rightarrow 0
\end{equation}

\noindent Since the cell structure of $\tilde X$ can be chosen so 
that all $C_i(\tilde X)$ are free $\pi_1(X)$-modules and 
since the structure of $\CC[\ZZ]$-module on 
$\CC[\pi_1(X)] \otimes_{\pi_1'} V_{\pi_1'}$ yields a 
free $\CC[\ZZ]$ module this implies that
$H_2(\tilde X_{\pi_1'},V_{\pi_1'}) 
\subseteq C_2(\tilde X) \otimes_{\pi_1'} V_{\pi_1'}$
is also free. Hence $\Delta_2(t)=1$.
Moreover, $\Delta_0(t)$ can be calculated in terms of $\rho$ as follows.
Denoting by $\partial_i: C_i(\tilde X) \rightarrow C_{i-1}(\tilde X)$ 
the boundary operator of the chain complex of the universal cover,
one has $H_0( \tilde X) \otimes_{\pi_1'} V_{\pi_1'}={\rm Coker} \partial_1 
\otimes {\rm id}$, $C_0( \tilde X)=\CC[\pi_1(X)]$, and 
${\rm Im} \partial_1$ is the augmentation ideal of the 
group ring $\CC[\pi_1(X)]$. Hence 
$H_0(\tilde X_{\pi_1'}) \otimes_{\pi_1'} V_{\pi_1'}
=\CC \otimes_{\pi_1'(X)} V_{\pi_1'}$
where the action of $\pi_1'(X)$ on $\CC$ is 
trivial. This implies that the latter is isomorphic 
to $\pi_1'(X)$ covariants of $V$. The $\CC[\ZZ]$-module structure
is given by the action of the generator of $\pi_1/\pi_1'$ on the 
$\pi_1'$ covariants $V^{\pi_1'}$ 
of $V$. For example, if $\rho$ has an abelian
image, i.e. $\pi_1' \subseteq \Ker \rho$, then $V^{\pi_1'}=V$ and 
$\Delta_0$ is the characteristic polynomial of the generator 
of $\ZZ$ acting on $V$ (cf. \cite{KL}).

Returning back to the case of an arbitrary CW-complex $X$, 
the relation between 
the twisted Alexander polynomial and twisted characteristic varieties 
can be stated as follows:

\begin{prop} Let $X$ be a CW complex such that $H_1(X,\ZZ)=\ZZ$,
$\rho: \pi_1(X) \rightarrow U(V)$ is a unitary local system 
and $\Ch_n^l \subseteq \Spec \CC[\ZZ]=\CC^*$ is the collection 
of characteristic varieties associated with $(X,\rho)$ in definition 
\ref{twistedcharvar}.
For $\xi \in \CC^*$ let $l_n(\xi)={\rm max} \{l \ \vert \xi \in \Ch^l_n(X,\rho)
\}$ and $b_n={\rm min} \{ l_n(\xi) \vert \xi \in \Spec \CC[\ZZ] \}$.
Then $\Delta_n(\xi)=0$ if and only if $l_n(\xi)>b_n$.
\end{prop}

\begin{proof} Let us consider the cyclic decomposition of 
$R=\CC[t,t^{-1}]$-module $H_n(X_{ab},V_{ab})$:

\begin{equation}\label{cyclicdecomp}
                H_n(X_{ab},V_{ab})=R^{\oplus c_n}
 \oplus R/(f_n^s) \ \ \
\end{equation}
Then $b_n$ is equal to 
the rank $c_n$. Morever 
for $\xi \in \Spec R=\CC^*$  
the integer $l_n(\xi)$ is equal to $c_n$ plus 
the number of torsion summands 
in (\ref{cyclicdecomp}) having
$\xi$ as its root. This yields the claim.
\end{proof}

\begin{prop}\label{proptwistedcharvar}
 If $\pi_i(X)=0$ for $2 \le i <k$ or if $k=1$
then 
the dimension of the homology group $H_k(X,V \otimes L_{\chi})$ corresponding 
to a character $\chi \in \Char\pi_1(X)$ and a local system $\rho$
is given by:
$$ {\rm max} \{ i \vert \chi \in \Ch_k^i(X,\rho) \}$$
\end{prop}

\begin{proof} Using the identification $C_*(\tilde X) \otimes_{\pi}
(V \otimes_{\pi} L_{\chi})=(C_*(\tilde X) \otimes_{\pi_1'} V) 
\otimes_{H_1(X)} \otimes L_{\chi}$ (recall that $C_*(\tilde X) 
\otimes_{\pi_1'} V$ has the structure of $\CC[H_1(X,\ZZ)]$-module:
$l(c \otimes v)=(c\bar l^{-1} \otimes \bar lv)$ where $l \in H_1(X,\ZZ)$
and $\bar l$ is its lift to $\pi_1$), the homology $H_*(X,V \otimes L_{\chi})$
can be obtained as the abutment of the spectral sequence:

\begin{equation}
H_p(H_1(X,\ZZ),H_q(C_*(\tilde X_ab) \otimes_{\pi'} V) \otimes_{\CC} L_{\chi})
 \Rightarrow H_{p+q}(X,V \otimes L_{\chi})
\end{equation}

Considering the lower degree terms we have:

\begin{equation}
H_{k+1}(H_1(X,\ZZ),L_{\chi}) \rightarrow
H_k(\tilde X_{ab},V)_{\chi, H_1(X,\ZZ)} \rightarrow  
H_k(X,V \otimes_{\pi_1} L_{\chi} 
\rightarrow H_k(H_1(X,\ZZ),L_{\chi})
\end{equation}

Since for $\chi \ne 1$ one has $H_k(H_1(X,\ZZ),\chi)=0$. 
Using the interpretion of the dimension of the middle term
via Fitting ideals 
we obtain the result. 

\end{proof}

\begin{theorem}\label{roots}
Let $X=\PP^2-C-L$ where $C$ is an irreducible curve
and $L$ is a line at infinity. Let $\rho$ be a unitary 
representation 
of the fundamental group and let $F$ be the extension of $\QQ$ 
generated by the eigenvalues of $\rho(\gamma)$ where $\gamma$
is a boundary of a small disk transversal to $C$ 
at its non singular point. Then the roots of $\Delta_{\rho}(C)$ 
belong to a cyclotomic extension of $F$.   
\end{theorem}

{\it Proof.} Let $f=0$ be the equation of the curve and let 
$\omega= A {{df} \over f}$ be the matrix of flat (logarithmic 
for a log resolution of the pair $(\PP^2,C \cup L)$) connection 
i.e. in a neighbourhood of any point in $\PP^2$, except for 
 those which are singularities of $C$, the connection 
is given by $\nabla_{\chi}(\cdot)=d \cdot + \omega \wedge \cdot$. 
Assume first that for a loop $\gamma$ having the linking number with $C$ equal
to one we have:  
\begin{equation}\label{unitary}
\chi(\gamma)=exp \ 2 \pi \sqrt{-1} \alpha
\end{equation}
where $\alpha \in \RR$ i.e. $\chi$ is unitary.
If $\pi: (X',E) \rightarrow (\PP^2,C\cup L)$ is the log resolution
(i.e., $E$ contains the exceptional set and the proper 
preimages of $C$ and $L$)
and  $ m_i={\rm mult}_{E_i} \pi^*(f)$ then the residue of 
log extension of flat connection on $V \otimes L_{\chi}$ 
is equal to $m_i(A(0)+\alpha I)$.
Hence the jumping loci are 

\begin{equation}\label{untiaryinequality}
m_i(\xi_{i,j}+\alpha)=n_i
\end{equation}
 for some 
integers $n_i$. Hence 
$(e^{2 \pi \sqrt{-1} \alpha})\zeta_{m_i}=e^{2 \pi \sqrt {-1} \xi_{i,j}}$ 
where $\zeta_{m_i}^{m_i}=1$
and the claim follows.

If the character $\chi$ is not unitary, i.e. ${\rm Im}\alpha \ne 0$, 
can one see that $dim H^1(V \otimes L_{\chi})=0$ as follows. First 
notice that the Deligne's Hodge deRham spectral sequence 
(\ref{hodgespectral}) has in  $E_1$ all terms equal to zero 
except possibly for the following: 
\begin{equation}
\begin{array}{ccc}
H^2(X',\V_{\chi}) & 0 
 & 0 \\
H^1(X',\V_{\chi}) & H^1(X',\Omega^1(\log E) \otimes \V_{\chi}) & 
0 \\
H^0(X',\V_{\chi}) & H^0(X',\Omega^1(\log E) \otimes \V_{\chi}) & 
H^0(X',\Omega^2(\log E) \otimes \V_{\chi})
\end{array}
\end{equation}
Indeed, the terms above the diagonal $i+j=2$ are zeros which 
can be seen as follows.
This spectral sequence, for variable $\chi$,
depends only on the locally trivial 
bundle supporting 
the Deligne's extension 
of the connection corresponding to local system $V \otimes L_{\chi}$.
This extension in turn does not depend on imaginary part 
${\rm Im} \alpha$ (but differential 
a priori may be dependent on the latter). Since in the case $\alpha 
\in \RR$ these terms are trivial due to degeneration 
of this spectral sequence in term $E_1$ and vanishing of
$H^i(X,V \otimes L_{\chi})$ for $i>2$ the claim follows. 
This spectral sequence
degenerates in term $E_1$ even if $\chi$ is not unitary since the differential
$d_1^{\alpha}$ in the spectral sequence corresponding to 
the character $\chi$ given by (\ref{unitary})
is the map induced on cohomology by the 
connection $\nabla_{\chi}: \Omega^i(\log E) \otimes \V 
\rightarrow \Omega^{i+1}(\log E) \otimes \V$ 
and hence depends on $\alpha$ linearly. Since it is trivial for 
$\alpha \in \RR$ (i.e. in the unitary case) 
it is trivial for all $\alpha$. The transgression
$d_2^{\alpha}$ depends on $\alpha$ linearly also  
as can be seen from the following description of $d_2^{\alpha}$: 
since for a form $\eta \in H^1(X',\V_{\chi})$ 
the differential  $d_1^{\alpha}$ takes its class to zero and hence 
we have $\nabla_{\chi}(\eta)=d\bar \eta$ where $\bar \eta$ is a $0$-form 
in $\Omega^1(\log E) \otimes \V_{\chi}$ and $\nabla_{\chi}(\bar \eta)$
represents the transgression.

\begin{remark} Hodge-deRham spectral sequence may not 
degenerate in term $E_1$ for non-unitary local system
(cf. \cite{Ramanan}).
\end{remark}

\begin{remark} If $\dim V=1$ (or, more generally, if $Im(\rho)$ is 
abelian), the twisted characteristic varieties coincide with the 
ordinary ones. If $\xi=\rho(\gamma)$ then the roots of $\Delta_{\rho}(t)$
are $\xi^{-1}\epsilon$ where $\epsilon$ is a root of ordinary 
Alexander polynomial (cf. also \cite{KL})
and the claim of \ref{roots} is immediate.
\end{remark}

One way to obtain a local system of higher rank is the following.
Let $X$ be a CW -complex and 
let $H$ be a normal subgroup of $\pi=\pi_1(X)$ such that $\rk H/H' <\infty$
($H'$ is the commutator subgroup of $H$).
The action of $\pi$ on $H$ by conjugation preserves $H'$ and hence 
we obtain an action on $H/H'$. Since for such action of $\pi$ 
the subgroup $H'$ acts trivially we also have the action of $\pi/H$.
We shall denote by 
$\rho_H$ the corresponding representation of $\pi$ (or $\pi/H$)
on $V_H=(H/H') \otimes \CC$. 

\begin{theorem}\label{restriction} Let $H$ be a finitely generated 
normal subgroup of 
$\pi=\pi_1(\CC^2-C)$ and $\gamma \in \pi$ an element in the conjugacy class 
of a loop which is the boundary of a small disk transversal to 
$C$.
The eigenvalues of the action of $\gamma$ 
on $(\pi' \cap H/(\pi' \cap H)') \otimes \CC$ belong to a cyclotomic
extension of the field $F$ generated by the eigenvalues of 
$\rho_H(\gamma)$ acting on $(H/H') \otimes \CC$.
\end{theorem}

\begin{proof} Recall that if $G$ is the covering group of a Galois 
cover $f: Y \rightarrow X$ and $V$ is a local system on $X$ then 
one has the Leray spectral sequence 
\begin{equation}H_p(G,H_q(Y,f^*(V))) \Rightarrow 
H_{p+q}(X,V)
\end{equation} 
We shall consider the exact sequence of low degree 
terms corresponding to the spectral sequence of $G=\pi'/\pi' \cap H$ 
acting on the 
covering space $X_{\pi' \cap H}$ of $X_{\pi'}$
corresponding to the subgroup $\pi' \cap H \subset \pi'$:
\begin{equation}
 H_p(G,H_q(X_{\pi' \cap H},f^*(V_H))) \Rightarrow H_{p+q}(X_{\pi},V_{H})
\end{equation} 
After the identification $H_1(X_{\pi' \cap H},f^*(V_H))=H_1(X_{\pi' \cap H})
\otimes_{\CC} V_H$ and $H_0(X_{\pi' \cap H})=\CC$
(since the pullback of the local system $V_H$ 
on $X_{\pi' \cap H}$ is trivial because the restriction of $\rho_H$ 
on $\pi'\cap H$ is trivial) we have:

\begin{equation}
 H_2(G,\CC) \rightarrow (H_1(X_{\pi' \cap H}) \otimes_{\CC} V_H)_G \rightarrow 
   H_1(X_{\pi'},V_H) \rightarrow H_1(G,\CC)
\end{equation}

The eigenvalues of $\rho(\gamma)$ 
belong to the cyclotomic extension of $F$ and since the 
action of $G$ is trivial on $H_1(X_{\pi' \cap H}))$ we obtain the claim.
\end{proof}

\begin{remark} In the case when $\pi/H$ is finite this result is not new:
if so, then $H$ is the fundamental group of a quasiprojective variety and 
the Alexander polynomial is cyclotomic in this case (\cite{ample}).
On the other hand this theorem imposes restrictions on finitely 
generated subgroups of the fundamental group
and this appears to be a non-trivial
restriction on this class of groups for which there 
exist a plane curve $C$ such that $\pi=\pi_1(\PP^2-C)$.
\end{remark}

\section{Examples}

\bigskip

\subsection{Points in $\PP^1$} Let us consider the case when $X$ is  
the complement to $r+1$ distinct points in $\PP^1$. Then $\pi_1(X)$ 
the free group on $r$ generators and $\Char\pi_1(X)={\CC^*}^r$.
Let $(\alpha_1,....,\alpha_r) \in \U$ (cf.(\ref{funddomain}))
and not all $\alpha_i=0$. 
Let $[-(\sum_i \alpha_i)]=-k$ 
($1 \le k \le r$). The Deligne's extension then is $\O_{\PP^1}(-k)$.
We have: 

\begin{equation}
\dim H^0(\Omega^1(\log \ D)(-k))=\dim H^0(\O(-2+r+1-k))=r-k
\end{equation}

\begin{equation}
\dim H^1(\O(-k))=k-1 
\end{equation}

Hence the dimensions of the graded components of Hodge filtration 
take constant values within the polytopes: 
\begin{equation}
    k < \sum \alpha_i \le k+1
\end{equation}
with $\dim H^1(\L_{\chi})=r-1$.

\subsection{Generic arrangement in $\PP^2$} The cohomology of a non-trivial 
local system on the complement $X$ to generic arrangement of $r+1$ lines 
are given 
as follows:

$$\rk H^2(X,\L_{\chi})= 1+{{r^2-3r} \over 2} \  \  \  
\rk H^i(X,\L_{\chi})=0\ \ \ 
 (i\ne 2) $$  

Together with the degeneration of the Hodge spectral sequence 
(\ref{hodgespectral})
this yields that each of the bundles of this spectral sequence 
in this case has only one non-vanishing cohomology group i.e.

\begin{equation}
           H^p(\Omega^q(\log D) \otimes \V)=(-1)^p
e(\Omega^q(\log D) \otimes \V)
\end{equation}
where $e(\A)$ denotes the holomorphic euler characteristic 
of a bundle $\A$.
Let us consider a local system for which the residues along r lines 
are $\alpha_1, ...., \alpha_r$ and 
\begin{equation}\label{polytopegeneric}
k-1 < \sum \alpha_i \le k
\end{equation}

Then the Deligne's extension is $\O(-k)$ and we have:

\begin{equation}
    e(\PP^2,\O(-k))={{(k-1)(k-2)}\over 2} 
\end{equation}

\begin{equation}
e(\PP^2, \Omega^2(\log D)(-k))=e(\PP^2,\O(r-2-k))=
{{(r-k-1)(r-k) \over 2}}
\end{equation}

For the calculation of $e(\PP^2,\Omega^1(\log D)(-k))$ 
let us use the Riemann-Roch.
The logarithmic bundle has the following Chern polynomial
(where $h \in H^2(\PP^2,\ZZ)$ is the generator and 
$c(\A)=\sum c_i(\A)t^i$ and $c_i(\A)$ are the Chern classes):
\begin{equation}
c(\Omega^1(\log D))=(1-ht)^{-r+2}
\end{equation}
(cf. \cite{DK}).
Therefore we have the following expression for the Chern character:
\begin{equation}
   ch(\Omega^2(\log D)(-k)=2+(r-2-2k)h+(k^2-(r-2)k-{{r-2}\over 2})h^2
\end{equation}
Hence, since the Todd class $td(\PP^2)=1+{3 \over 2}h+h^2$, the 
 Riemann Roch yields:
\begin{equation}
e(\PP^2,\Omega^1(\log D)(-k))=k^2+r-k-rk
\end{equation}

Hence we obtain the following:

\begin{theorem}\label{generic}
 For the generic arrangement in $\PP^2$ one has the 
following:
$$\dim Gr_F^0={{(k-1)(k-2)}\over 2}, \ \ \dim Gr_F^1=rk+k-r-k^2, \ \ 
\dim Gr_F^2={{(r-k-1)(r-k) \over 2}}
$$
\end{theorem}

\subsection{Cone over generic arrangements} 
In the case of the cone over $r$ points in $\PP^1$ we obtain
an ordinary plane curve singularity 
of multiplicity $r$ which is the case discussed already in \cite{alexhodge}.
Consider now the cone over the generic arrangement in $\PP^2$ considered 
in the previous example, i.e. the 
arrangement $C\A$ of $r+2$ planes in $\PP^3$ with $r+1$ planes 
$H_1,...,H_{r+1}$
forming 
an isolated non-normal crossing and the remaining plane $H_{\infty}$ 
(at infinity)
being transversal to the first $r+1$ planes. One can make calculations on 
the blow up of $\PP^3$ at the non-normal crossing point of the arrangement
but from 
$$0 \rightarrow \ZZ \rightarrow \pi_1(\PP^3-C\A) \rightarrow 
\pi_1(\PP^2-\A) \rightarrow 0 $$ 
one has the identification of $\Char \pi_1(\PP^2-\A)$ with 
the subgroup of $\Char\pi_1(\PP^3-C\A)$ in the standard coordinates in the 
latter given by $$x_1 \cdot ... \cdot x_{r+1}=1$$
From the Gysin sequence associated with 
$H_{\infty}-\cup_{i=1}^{i=r} H_i \cap H_{\infty} \subset \PP^3-
\cup_{i=1}^{i=r} H_i$
 we have the identification of the 
Hodge structures :
\begin{equation}
    H^3(\PP^3-C\A, \L_{\chi})=H^2(H_{\infty}-H_{\infty} \cap \cup H_i,
\L_{\chi})
\end{equation}
 
Hence the polytopes $\S^2_k$ are given by 
$\alpha_1+...+\alpha_{r+1}=k$ and in each of these polytopes 
we have the value of $\dim Gr_F^k$ given by the theorem \ref{generic}.

\subsection{Arrangement from a pencil of plane quadrics}
Consider the arrangement of six lines $L_1,..,L_6$ composed of the three
sides 
of triangle and the three cevians in it. Let $H_{\infty}$ be the 
line at infinity. $\Char \pi_1(\PP^2-\cup L_i \cup H_{\infty})={\CC^*}^6$.
The characteristic variety is given by
\begin{equation}\label{cevachar}
x_{i_1}\cdot x_{i_2} \cdot x_{i_3}=1
\end{equation}
where $(i_1,i_2,i_3)$ runs through triples having a triple point in common
(cf. \cite{abcov}).
This is a two dimensional sub-torus in the torus of characters.
Let us consider the characters $exp(u)$ where $u=(\alpha_1,...,\alpha_6)$
such that $\sum \alpha_{i_1}+\alpha_{i_2}+\alpha_{i_3}=p+1$ ($p=0,1$)
$(i_1,i_2,i_3)$ as in (\ref{cevachar}). Let $\tilde \PP^2$ be the blow 
up of $\PP^2$ at the vertices of the triangle. 
The Deligne's
extension is $\O(\sum_i(p+1)E_i-2pH_{\infty})$.
We have using $K_{\tilde \PP^2}=\sum E_i-3H_{\infty}$ and  Serre duality
$$H^1(\tilde \PP^2,\O(\sum_i 2E_i-4H_{\infty}))=H^1(\tilde \PP^2,\O(-\sum E_i+
H_{\infty}))=H^1(\PP^2,\J_B(H_{\infty}))=1$$
where $\J_B$ is the ideal sheaf of the set of the vertices of the triangle
(the last equality follows from the Cayley-Bacharach theorem, 
cf. \cite{abcov}). Similarly one obtains ($\V$ is the Deligne's extension):
 $$Gr_F^0={\rm dim} H^0(\tilde \PP^2, \Omega^1({\rm log} \cup E_i 
\cup H_{\infty})\otimes \V)=1$$

\subsection{Arrangement from a net of quadrics in $\PP^3$}
Consider the arrangement $D_{8,4}=\cup_{i=1,...8} H_i$ 
of eight planes in $\PP^3$ introduced in \cite{ample}.
Recall that these eight planes split in four pairs forming $4$ quadrics
belonging to a net. This net has eight base points which are 
the  {\it only} eight  non-normal 
crossings in the arrangement. Let $H_{\infty}$ be a generic hyperplane in 
$\PP^3$. One has $\Char \pi_1(\PP^3-D_{8,4} \cup H_{\infty})={\CC^*}^8$.
Consider the subset in $\U$ (cf. section \ref{quasiproj}) given 
by 
\begin{equation}\label{equationofhodge}
\bar S^2_3: \ \ \ \alpha_{i_1}+...+\alpha_{i_4}=3 
\end{equation}
where $(i_1,...,i_4)$
 runs through eight unordered collections of indices such that 
$H_{i_1},...,H_{i_4}$ form a quadruple of planes containing 
one of eight non-normal crossings points. Note that the dimension (over $\RR$)
of the set of solutions is equal to three and the Zariski closure of 
this set is support of the homotopy group 
${\rm Supp} \pi_2(\PP^3-D_{8,4}) \otimes \CC$ 
considered in \cite{ample}. We claim that for any 
$u=(\alpha_1,...,\alpha_8) \in \bar S^2_3 $ and for the 
corresponding character $exp(u)$ one has 
\begin{equation}\label{lastclaim}
\dim Gr^3_FH^2(\PP^3-D_{8,4} \cup H_{\infty},\L_{exp(u)})=1
\end{equation}
Indeed let $\tilde \PP^3$ be the blow up of $\PP^3$ in eight base points
and $E_1,...,E_8$ the corresponding exceptional components.
For a character $exp(u)$, where $u$ belongs to (\ref{equationofhodge}),
the Deligne's extension is $\O(3E_1+...+3E_8-6H_{\infty})$, since
each of the eight planes in the arrangement contains four base 
points. The canonical class of $\tilde \PP^3$ is $2E_1+...+2E_8-4H_{\infty}$ 
and therefore we have:
\begin{equation}
\dim Gr^3_FH^2(\PP^3-D_{8,4} \cup H_{\infty},\L_{exp(u)})=
\dim H^2(\O((3E_1+...+3E_8-6H_{\infty})=
\end{equation}
(using Serre duality)
$$\dim H^1(\O(-E_1-...-E_8+2H_{\infty})=\dim H^1(\PP^3,\J_B(2H))$$
where $\J_B$ is the ideal sheaf of the base locus. Since $B$ is the complete 
intersection of three quadrics in $\PP^3$, the claim (\ref{lastclaim})
follows from the 
Cayley-Bacharach theorem. 
In fact $\dim Gr_F^p=1$ for $p=0,1$ as well and the corresponding 
characters are the exponents of $u$ for which 
$\alpha_{i_1}=...=\alpha_{i_r}=p+1$.

\end{document}